\documentclass[leqno]{article}
\usepackage{latexsym,amsmath,amsfonts,mathrsfs,amssymb,amsthm}
\usepackage[usenames]{color}
\usepackage{graphicx}
\def\R{\mathbb R}
\def\N{\mathbb N}

\def\Z{\mathbb Z}

\def\a{\alpha}
\def\e{\epsilon}
\def\d{\delta}
\def\Y{\mathbb Y}
\def\T{\mathbb T}
\def\P{\mathbb P}
\def\be{\begin{equation}}
\def\ee{\end{equation}}
\def\bs{\backslash}
\def\qed{\hfill$\Box$\bigskip}
\def\nd{\noindent Proof. }
\numberwithin{equation}{section}
\newtheorem{lem}[equation]{Lemma}
\newtheorem{pro}[equation]{Proposition}
\newtheorem{defn}[equation]{Definition}

\newtheorem{cor}[equation]{Corollary}
\newtheorem{rem}[equation]{Remark}
\begin{document}
\bigskip

\centerline{\Large \textbf{Global regularity for minimal sets near a $\T$ set and counterexamples}}

\bigskip

\centerline{\large Xiangyu Liang}

\vskip 1cm

\centerline {\large\textbf{Abstract.}}

\bigskip

We discuss the global regularity for 2 dimensional minimal sets that are near a $\T$ set, that is, whether every global minimal set in $\R^n$ that looks like a $\T$ set at infinity is a $\T$ set or not. The main point is to use the topological properties of a minimal set at large scale to control its topology at smaller scales. This is the idea to prove that all 1-dimensional Almgren-minimal sets in $\R^n$, and all 2-dimensional Mumford-Shah minimal sets in $\R^3$ are cones. In this article we discuss two types of 2-dimensional minimal sets: Almgren-minimal set in $\R^3$ whose blow-in limit is a $\T$ set; topological minimal sets in $\R^4$ whose blow-in limit is a $\T$ set. For the first one we eliminate an existing potential counterexample that was proposed by several people, and show that a real counterexample should have a more complicated topological structure; for the second we construct a potential example using a Klein bottle.
\bigskip


\bigskip


\setcounter{section}{-1}

\section{Introduction}

This paper deals with the global regularity of two-dimensional minimal sets that looks like a $\T$ set at infinity in $\R^3$ and $\R^4$. The motivation is that we want to decide whether all global minimal sets in $\R^n$ are cones. 

This Bernstein type of problems is a typical interest for all kinds of minimizing problems in geometric measure theory and calculus of variations. It is natural to ask how does a global minimizer look like, as soon as we know already the local regularity for minimizers. Well known examples are the global regularity for complete 2-dimensional minimal surfaces in $\R^3$, area or size minimizing currents in $\R^n$, or global minimizers for the Mumford-Shah functional. Some of them admit very good descriptions for global minimizers. See \cite{Be, Mo89, Mo86,DMS} for further information.



Now let us say something more precise about minimal sets. Briefly, a minimal set is a closed set which minimizes the Hausdorff measure among a certain class of competitors. Different choices of classes of competitors give different kinds of minimal sets. So we have the following general definition.

\begin{defn}[Minimal sets]\label{min}
Let $0<d<n$ be integers. A closed set $E$ in $\R^n$ is said to be minimal of dimension $d$ in $\R^n$ if 
\be H^d(E\cap B)<\infty\mbox{ for every compact ball }B\subset \R^n,\ee
and
\be H^d(E\bs F)\le H^d(F\bs E)\ee
for any competitor $F$ for $E$.
\end{defn}

\begin{rem}We can of course give the definition of locally minimal sets, where we merely replace $\R^n$ in Definition 0.1 (and the definitions of Almgren and topological competitors which will appear later) by any open set $U\subset\R^n$. This makes no difference when we discuss local regularity. But for global regularity, the ambient space $\R^n$ always plays an important role.
\end{rem}

In this paper we will discuss the following two kinds of minimal sets, that is, sets that minimize Hausdorff measure among two classes of competitors.

\begin{defn}[Almgren competitor (Al-competitor for short)] Let $E$ be a closed set in $\R^n$. An Almgren competitor for $E$ is a closed set $F\subset \R^n$ that can be written as $F=\varphi(E)$, where $\varphi:\R^n\to\R^n$ is a Lipschitz map such that there exists a compact ball $B\subset\R^n$ such that
\be \varphi|_{B^C}=id\mbox{ and }\varphi(B)\subset B.\ee 
Such a $\varphi$ is called a deformation in $B$, and $F$ is also called a deformation of $E$ in $B$.
\end{defn}

Roughly speaking, we say that $E$ is Almgren-minimal when there is no deformation $F=\varphi(E)$, where $\varphi$ is Lipschitz and $\varphi(x)-x$ is compactly supported, for which the Hausdorff measure $H^d(F)$ is smaller than $H^d(E)$. The definition of Almgren minimal sets was invented by Almgren \cite{Al76} to describe the behaviors of physical objects that span a given boundary with as little surface area as possible, such as soap films.  

The second type of competitors was introduced by the author in \cite{topo}, where she tried to generalize the definition of Mumford-Shah minimal sets (MS-minimal for short) to higher codimensions. In both definitions for MS competitors and topological competitors, we ask that a competitor keeps certain topological properties of the initial set. Sometimes this condition is easier to handle than the deformation condition that is imposed for Al-competitors.

\begin{defn}[Topological competitor] Let $E$ be a closed set in $\R^n$. We say that a closed set $F$ is a topological competitor of dimension $d$ ($d<n$) of $E$, if there exists a ball $B\subset\R^n$ such that

1) $F\bs B=E\bs B$;

2) For every Euclidean $n-d-1$-sphere $S\subset\R^n\bs(B\cup E)$, if $S$ represents a non-zero element in the singular homology group $H_{n-d-1}(\R^n\bs E;\Z)$, then it is also non-zero in $H_{n-d-1}(\R^n\bs F;\Z)$.
\end{defn}

\begin{rem}When $d=n-1$, this is the definition of MS-competitor, where we impose a separation condition for the complementary of the set.
\end{rem}

The therefore defined class of topological minimizers is contained in the class of Almgren minimal sets (c.f.\cite{topo}, Corollary 3.17), and admits some good properties that we are not able to prove for Almgren minimal sets.

\bigskip

Our goal is to show that a minimal set in $\R^n$ is a cone. Topological minimal sets are automatically Almgren minimal, hence almost all the information given below is for Almgren minimal sets.

Let $E$ be a $d-$dimensional reduced Almgren minimal set in $\R^n$. Reduced means that there is no unnecessary points. More precisely, we say that $E$ is reduced when
\be H^d(E\cap B(x,r))>0\mbox{ for } x\in E\mbox{ and }r>0.\ee

Recall that the definition of minimal sets is invariant modulo sets of measure zero, and it is not hard to see that for each Almgren (resp. topological) minimal set $E$, its closed support $E^*$ (the reduced  set $E^*\subset E$ with $H^2(E\bs E^*)=0$) is a reduced Almgren (resp. topological) minimal set. Hence we can restrict ourselves to discuss only reduced minimal sets.

Now fix any $x\in E$, and set
\be \theta_x(r)=r^{-d}H^d(E\cap B(x,r)).\ee

Then this density function $\theta_x$ is nondecreasing for $r\in ]0,\infty[$ (c.f.\cite{DJT} Proposition 5.16). In particular the two values
\be \theta(x)=\lim_{t\to 0^+}\theta_x(t)\mbox{ and }\theta_\infty(x)=\lim_{t\to \infty}\theta_x(t)\ee
exist, and are called density of $E$ at $x$, and density of $E$ at infinity respectively. Notice that $\theta_\infty(x)$ does not depend on $x$, hence we denote it by $\theta_\infty$.

Now Theorem 6.2 of \cite{DJT} says that if $E$ is a minimal set, $x\in E$, and $\theta_x(r)$ is constant on $r$, then $E$ is a minimal cone centered on $x$. Thus by the monotonicity of the density functions $\theta_x(r)$ for any $x\in E$, if we can find a point $x\in E$ such that $\theta(x)=\theta_\infty$, then we are done.

On the other hand, the possible values for $\theta(x)$ and $\theta_\infty$ for any $E$ and $x\in E$ are not arbitrary. 
By Proposition 7.31 of \cite{DJT}, for each $x$, $\theta(x)$ is equal to the density at the origin of a $d-$dimensional Al-minimal cone in $\R^n$. An argument around (18.33) of \cite{DJT}, which is similar to the proof of Proposition 7.31 of \cite{DJT}, gives that $\theta(x)$ is also equal to the density at the origin of a $d-$dimensional Al-minimal cone in $\R^n$. In other words, if we denote by $\Theta_{d,n}$ the set of all possible numbers that could be the density at the origin of a $d-$dimensional Almgren-minimal cone in $\R^n$, then $\theta_\infty\in\Theta_{d,n}$, and for any $x\in E$, $\theta(x)\in\Theta_{d,n}$.

Thus we restrict the range of $\theta_\infty$ and $\theta(x)$. Recall that the set $\Theta_{d,n}$ is possibly very small for any $d$ and $n$. For example, $\Theta_{2,3}$ contains only three values: 1 (the density of a plane), 1.5 (the density of a $\Y$ set, which is the union of three closed half planes with a common boundary $L$, and that meet along the line $L$ with $120^\circ$ angles), and $d_T$ (is the density of a $\T$ set, i.e., the cone over the 1-skeleton of a regular tetrahedron centered at 0). (See the figure below).

\centerline{\includegraphics[width=0.32\textwidth]{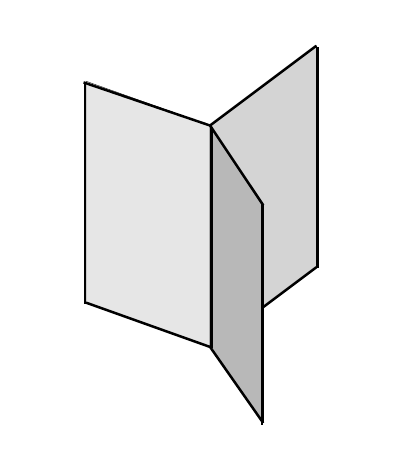} \includegraphics[width=0.4\textwidth]{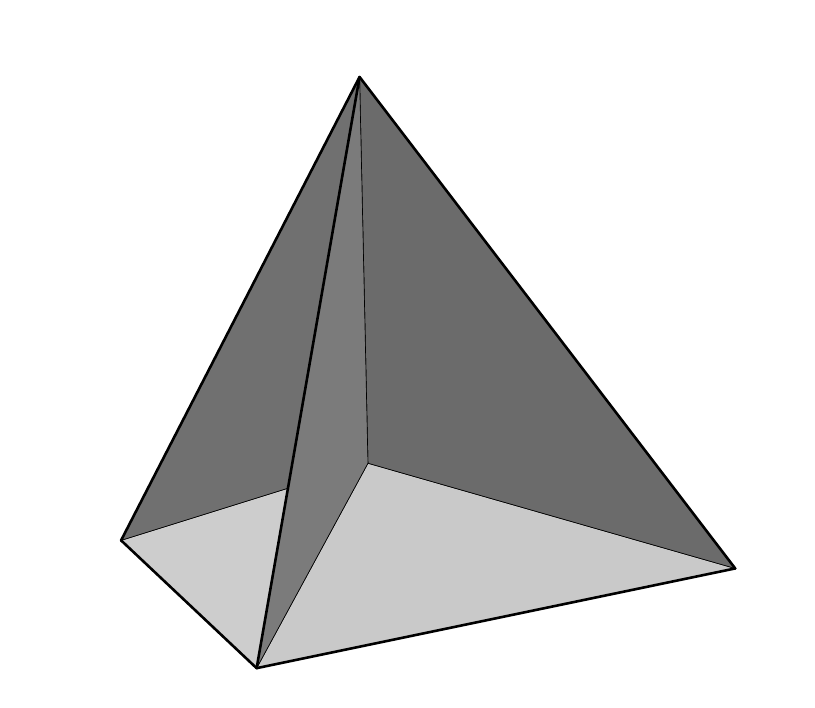}}

Recall that the reason why $\theta_\infty$ has to be in $\Theta_{d,n}$ is that, for any Al-minimal set $E$, all its blow-in limits have to be Al-minimal cones (c.f. Argument around (18.33) of \cite{DJT}). A blow-in limit of $E$ is the limit of any converging (for the Hausdorff distance) subsequence of
\be E_r=r^{-1}E, r\to\infty .\ee

Hence the value of $\theta_\infty$ implies that at sufficiently large scales, $E$ looks like an Al-minimal cone of density $\theta_\infty$. 

This is the same reason why $\theta(x)\in\Theta_{d,n}$. Here we look at the behavior of $E_r$ when $r\to 0$, and the limit of any converging subsequence is called a blow-up limit (this might not be unique!). Such a limit is also an Al-minimal cone $C$ (c.f. \cite{DJT} Proposition 7.31). This means, at some very small scales around each $x$, $E$ looks like some Al-minimal cone $C$ of density $\theta(x)$. In this case we call $x$ a $C$ point of $E$.

After the discussion above, our problem will reduce to the following one: do all minimal cones $C$ satisfy the following property ?
\be \begin{split}&\mbox{If a }\mbox{minimal set }E\mbox{ looks very much like }C\mbox{ at some scale,}\\
&\mbox{ then }E\mbox{ contains a point }x\mbox{ such that the density }\theta(x)\mbox{ is the same as that of }C\mbox{ at the origin}.\end{split}\ee

The discuss above uses only the values of densities at small scale and at infinity. But we can interpret it geometrically in the following way:
the blow-up limit $C_x$ of $E$ at $x$ should be: either $C$ itself (which might be too nice to hope), either some minimal cone that admits the same density as $C$. The latter, at least in all known cases, is equivalent to say that $C_x$ should admit the same topology as $C$. Thus a geometrical restatement of (0.13) is
\be\begin{split}&\mbox{If a }\mbox{minimal set }E\mbox{ looks very much like }C\mbox{ at some scale,}\\
&\mbox{ then }E\mbox{ contains a point }x\mbox{ which has a blow-up limit that admits the same topology as }C.\end{split}\ee

This seems to be a topological property that depends on each minimal cone $C$. In fact, topological arguments work when $C$ is a 2-dimensional $\Y$ set in $\R^3$ (c.f. \cite{DJT} Proposition 16.24).  However for the $\T$ set, no topological argument is enough. There is an example $E_0\subset\overline B(0,1)$ proposed by several persons (c.f. \cite{Mo95} page 110, or \cite{DJT} section 19), which satisfies that
\be E_0\cap\partial B(0,1)=T\cap\partial B(0,1),\ee
where $T$ is a $\T$ set centered at the origin, and $E_0$ satisfies all the known local regularity properties for Al-minimal sets, but $E_0$ contains no $\T$ point. (c.f.\cite{DJT}, Section 19 for a description for $E_0$).

We will talk about this later in this article. 

Besides the global regularity, the property (0.13) helps also to use the distance between a minimal set and a minimal cone to control the local speed of decay of the density function $\theta_x(r)$, because this property gives a lower bound of $\theta_x(r)$. When the speed is small, it might lead to nicer local regularity result. See \cite{DEpi} for more detail.

So in this paper we will treat the property (0.13) for the minimal cones $\T$ under the two types of definitions for "minimal".

In the first part of this paper, we discuss (0.13) for $\T$ sets in the senses of Almgren-minimal in $\R^3$. We eliminate the potential counterexample $E_0$ proposed above, and give some descriptions for real potential counterexamples if they exist. However we want to point out that this description asks that such a counterexample must admit a very complicated topology (c.f. Proposition 2.10 and Corollary 2.14), which, from the author's point of view, contradicts the spirit of minimal sets. In spite of that, we will still give an example of a set that admits this complicated topology. Hence we are still not able to prove (0.13) for $\T$ set in the sense of Almgren.

The discussion for the existence of a $\T$ point for topological minimal sets in $\R^4$  (since the problem in $\R^3$ has been already solved, c.f. \cite{DJT}) is done in Part II. The topological minimality seems to be stronger than Al-minimality, and it is proved in Section 18 of \cite{DJT} that in $\R^3$, the property (0.13) holds for $\T$ sets. In particular Theorem 1.9 in \cite{DJT} says that all 2-dimensional topological minimal sets in $\R^3$ are cones. However in $\R^4$, when the codimension is 2, things are complicated. We do not know the list of minimal cones in this case. And even if we make some additional assumption (See (5.1), which says that in $\R^4$ there is no 2-dimensional minimal cone whose density is less than $\T$ other than a plane or a $\Y$ set), we still end up with a topological counterexample that satisfies all the known local regularity properties. 

In this article, some of the results and arguments cited in \cite{DJT} exist also in some other (earlier) references, e.g. \cite{Ta}. But for simplify the article, the author will cite \cite{DJT} systematically throughout this article.

\textbf{Acknowledgement:} I would like to thank Guy David for many helpful discussions and for his continual encouragement. Part of the results in this paper were part of the author's doctoral dissertation at University of Paris-Sud 11 (Orsay). This work was partially supported by grants from R\'egion Ile-de-France.

\noindent\textbf{Some useful notation}

In all that follows, minimal set means Almgren minimal set;

$[a,b]$ is the line segment with end points $a$ and $b$;

$[a,b)$ is the half line with initial point $a$ and passing through $b$;

$B(x,r)$ is the open ball with radius $r$ and centered on $x$;

$\overline B(x,r)$ is the closed ball with radius $r$ and center $x$;

$\overrightarrow{ab}$ is the vector $b-a$;

$H^d$ is the Hausdorff measure of dimension $d$ ;

$d_H(E,F)=\max\{\sup\{d(y,F):y\in E,\sup\{d(y,E):y\in F\}\}$ is the Hausdorff distance between two sets $E$ and $F$.

$d_{x,r}$ : the relative distance with respect to the ball $B(x,r)$, is defined by
$$ d_{x,r}(E,F)=\frac 1r\max\{\sup\{d(y,F):y\in E\cap B(x,r)\},\sup\{d(y,E):y\in F\cap B(x,r)\}\}.$$

\part{Existence of a point of type $\T$ for 2-dimensional Al-minimal sets in $\R^3$}

\section{Introduction}

In this part we treat the old problem of the characterization of 2-dimensional Al-minimal sets in $\R^3$, and restrict the class of potential Al-minimal sets that are not cones.

Recall that this problem for 2-dimensional topological minimal sets in $\R^3$ (which coincides with MS-minimal sets in this case) has been positively solved in \cite{DJT}, where Theorem 1.9 says that all 2-dimensional MS-minimal sets in $\R^3$ are cones. The proof of this theorem is essentially to prove the property (0.13) for all 2-dimensional MS-minimal cones in $\R^3$. There are only three types of minimal cones in this case, which are the planes, the $\Y$ sets, and the $\T$ sets. In \cite{DJT}, (0.13) has been proved for planes and $\Y$ sets, only under the assumption of Almgren minimality. The MS-minimality is used to prove (0.13) for $\T$ sets, where Al-minimality seems to be less powerful.

But in Section 2 we are going to eliminate the well known potential counterexample (c.f. \cite{DJT} Section 19). Topologically this example satisfies all known local regularity properties for Al-minimal sets, but we will still manage to give another topological criterion (Proposition 2.10 and Corollary 2.14) for minimal sets that look like a $\T$ set at infinity, and use this property to prove that the potential counterexample, as well as some other similar sets, cannot be Almgren-minimal. This topological criterion seems to be really strange, and, intuitively, could not be satisfied by any global minimal sets. However topologically, sets that admit such a property exist, and we are going to construct such an example in Section 3.

In Section 4 we will treat another similar problem, that is, for a $\T$ set $T$, is the set $T\cap B(0,1)$ the only minimal set $E$ in $\overline B(0,1)$ such that $E\cap\partial B(0,1)=T\cap \partial B(0,1)$ ? While all the above arguments give some methods for controlling the measure of a set by topology, in Section 4 we will give some way to control the topology of a set by its measure.

\section{A topological criterion for potential counterexamples}

In this section we will give a topological necessary condition for potential 2-dimensional Almgren-minimal sets in $\R^3$ that are not cones.

First let us recall some known facts about such a set. Let $E$ be such a set. We look at the sets
\be E(r,x)=\frac 1r (E-x)\ee
where $r$ tends to infinity.

For every sequence $\{t_k\}_{k\in\N}$ which tends to infinity and such that $E(t_k,x)$ converges (in all compact sets, for the Hausdorff distance), the limit (called a blow-in limit) should be a minimal cone $C$ (c.f. \cite{DJT}, arguments around (18.33)). Now by the classification of singularities \cite{Ta}, $C$ should be a plane, a $\Y$ set or a $\T$ set. By \cite{DJT}, $C$ could not be a plane or a $\Y$ set. Hence $C$ is a $\T$ set. And thus there exists a $\T$ set $T$ centered at the origin, and a sequence $\{t_k\}_{k\in\N}$ such that
\be \lim_{k\to\infty} t_k=\infty\mbox{ and }\lim_{t_k\to\infty}d_{0,t_k}(E,T)=0.\ee

Denote by $B=B(0,1)$ the unit ball. Denote by $y_i,1\le i\le 4$ the 4 points of type $\Y$ of $T\cap\partial B$. Denote by $C$ the convex hull of $\{y_i,1\le i\le 4\}$, which is a regular tetrahedron inscribed in $B$. Set $T_C=T\cap C$. Then a simple calculation gives 
\be \frac 12H^2(\partial C)=\frac43\sqrt 3< 2\sqrt 2=H^2(T_C).\ee

Set $\delta=\frac 14(H^2(T_C)-\frac12 H^2(\partial C))$. Then a minor modification of the proof of Lemma 16.43 of \cite{DJT} gives
\begin{lem}There exists $\e_1>0$ such that if $d_{0,2}(E,T)<\e_1$, then
\be H^2(E\cap C)>H^2(T_C)-\d.\ee
\end{lem}

On the other hand, there exists $\e_2>0$ such that if $d_{0,1}(E,T)<\e_2$, then in the annulus $B(0,\frac 32)\bs B(0,\frac 12)$, $E$ is a $C^1$ version of $T$. (c.f. \cite{DJT} Section 18). More precisely, in $B(0,\frac 32)\bs B(0,\frac 12)$, the set $E_Y$ of points of type $\Y$ in $E$ is the union of four $C^1$ curves $\eta_i,1\le i\le 4$. Each $\eta_i$ is very near the half line $[o,y_i)$, and around each $\eta_i$, there exists a tubular neighborhood ${\cal T}_i$ of $\eta_i$, which contains $B([0,y_i),r)$ for some $r>0$, such that $E$ is a $C^1$ version of a $\Y$ set in ${\cal T}_i$. And for the part of $E\bs E_Y$, $E\cap B(0,\frac 32)\bs B(0,\frac 12)$ is composed of 6 flat surfaces $E_{ij},1\le i<j\le 4$. Each $E_{ij}$ is very near $T_{ij}$, where $T_{ij}$ is the cone over the great arc $l_{ij}$, which is the great arc on $\partial B$ that connects $y_i$ and $y_j$. Thus each $E_{ij}$ is a locally minimal set that is near a plane. Then by an argument similar to the proof of Proposition 6.14 of \cite{2p}, outside $\cup_{1\le i\le 4}{\cal T}_i$, $E_{ij}$ is the graph of a $C^1$ function of $T_{ij}$. Hence all in all, in $B(0,\frac 32)\bs B(0,\frac 12)$, $E$ is the image of $T$ by a $C^1$ diffeomorphism  $\varphi$, whose derivative is very near the identity.  

Thus by (2.2), and possibly modulo a dilation, we can suppose that for $t_k=2$,
\be d_{0,2}(E,T)<\min\{\e_1,\e_2\},\ee
which gives (2.5), and that in $B(0,\frac 32)\bs B(0,\frac 12)$, $E$ is a $C^1$ version of $T$.

In particular, on the boundary of $C$, $E\cap\partial C$ admits the same topology as $T\cap\partial C$. That is, $E\cap\partial C$ is composed of six piecewise $C^1$ curves $w_{ij},1\le i<j\le 4$, which meet at 4 endpoints $b_i,1\le i\le 4$, by sets of three. Each $b_i$ is very near $y_i$, where $y_i$ are the 4 points of type $\Y$ of $T\cap\partial C$. Denote by $w_{ij}$ the curve in $E\cap\partial C$ that connects $b_i$ and $b_j$. Then $w_{ij}$ is very near $[b_i,b_j]$. Moreover, if we denote by $\Omega_i,1\le i\le 4$, the connected component of $\partial C\bs E$ which is opposite to $b_i$, bounded by the $w_{kl},k,l\ne i$, then we can ask that $\e$ is small enough so that
\be \mbox{for each }1\le i\le 4, H^2(\Omega_i)>\frac14 H^2(\partial C)-\d,\ee
 where $\frac14 H^2(\partial C)$ is the measure of a face of $\partial C$. (See Figure 2-1.)

\centerline{\includegraphics[width=0.7\textwidth]{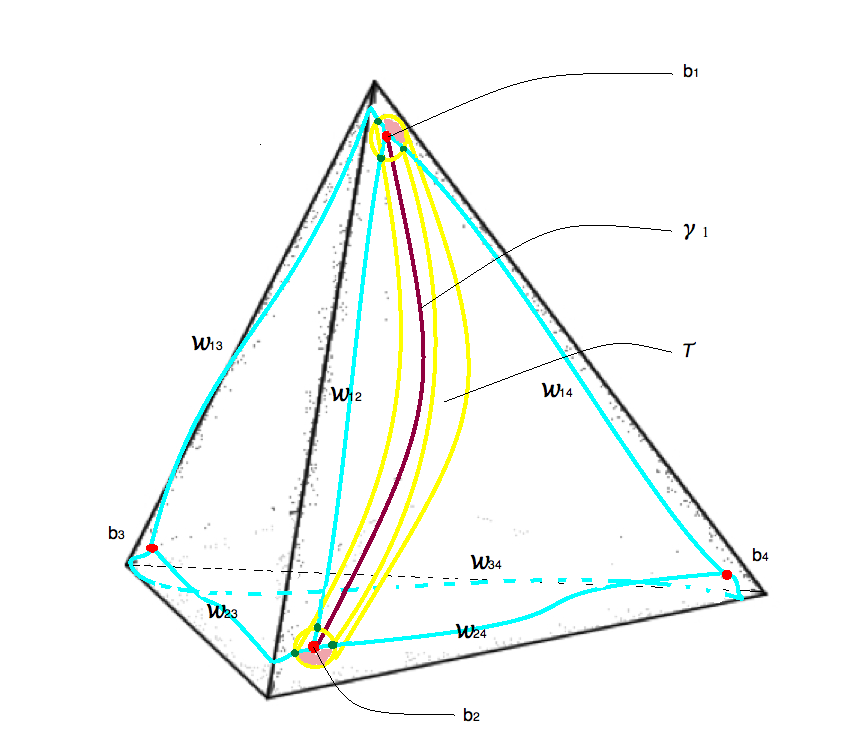}}
\nopagebreak[4]
\centerline{2-1}

Now suppose that there is no point of type $\T$ in $E\cap C$. Recall that $E_Y$ is the set of all points of type $\Y$ in $E$. Then by the $C^1$ regularity around points of type $\Y$ (c.f.\cite{DEpi} Theorem 1.15 and Lemma 14.6), $E_Y\cap C$ is composed of $C^1$ curves, whose endpoints are $b_i,1\le i\le 4$. Then there exists two curves $\gamma_1,\gamma_2\subset E_Y$ whose endpoints are the $b_i$. Suppose, for example, that $\gamma_1\cap\partial C=\{b_1,b_2\}$, and $\gamma_2\cap\partial D=\{b_3,b_4\}$.

Now by the $C^1$ regularity for points of type $\Y$ (\cite{DEpi} Theorem 1.15 and Lemma 14.6), for each $x\in\gamma_1$, there exists a neighborhood $B(x,r)$ such that in $B(x,r)$, $E$ is a $C^1$ version of $Y+x$, which cuts $B(x,r)$ into 3 connected components. Then by the compactness of $\gamma_1$, there exists $r>0$ such that in the tubular neighborhood $B(\gamma_1,r)$ of $\gamma_1$, $E$ is a distorted $Y$ set, whose singular set is $\gamma_1$, and $E$ divides $B(\gamma_1,r)$ into three connected components. Each component is a long tube that joins one of the three $\Omega_i$ near $b_1$ to one of the three $\Omega_i$ near $b_2$. Notice that if for $i\ne j$, $\Omega_i$ and $\Omega_j$ are connected by one of these long tubes, then they lie in the same connected component of $B\bs E$. As a result, there exist $1\le i,j\le 4, i\ne j$ such that $\Omega_i$ and $\Omega_j$ are in the same connected component of $B\bs E$, and there exists a long tube $\cal T$ along $\gamma_1$ which connects $\Omega_i$ and $\Omega_j$.

\smallskip

Now suppose that there exists a deformation $f$ in $C$ (see Definition 0.5), two indices $1\le i\ne j\le 4$, and two points $x\in \Omega_i,y\in\Omega_j$, such that
\be f(E)\subset C\bs[x,y].\ee

It is then not hard to find a Lipschitz deformation $g:C\bs B([x,y],r)\to G:=\partial C\bs (\Omega_i\cup\Omega_j)$ and $g=Id$ on $G$. For the construction of such a $g$, we can imagine that we enlarge the "hole" $B([x,y],r)$ and push every point in $C\bs B([x,y],r)$ towards the set $G$. For example we give, in Figure 2.2, a sketch for what happens, when $E\cap\partial C=T\cap\partial C$. For any set $E$ we have only to do some tiny modification, since $E\cap\partial C$ is a $C^1$ version of $T\cap\partial C$. Here for each half plane $D$ that is bounded by the line containing $[x,y]$, we just map $D\cap C\bs B([x,y],r)$ to the boundary $G\cap D$ (the thicker segments or point in the figure).

\centerline{\includegraphics[width=0.8\textwidth]{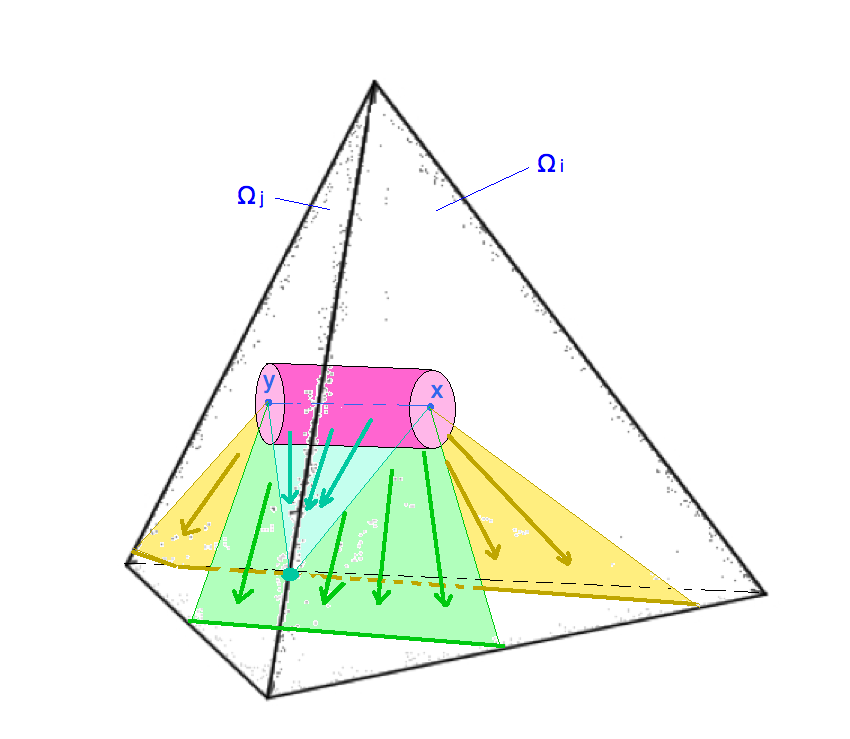}}
\nopagebreak[4]
\centerline{2-2}
 
 Then the function $h:=g\circ f$ sends $E$ to a subset of $G$ for $t$ large, and moreover, $g$ does not move $E\cap C=\cup_{k=1}^4\partial\Omega_k$. 

The above argument implies that in $C$ we can deform $E$ to a subset of $G$. 

Now by (2.5) and (2.7),
\be \begin{split}H^2(h(E))&\le H^2(G)=H^2(\partial C)-H^2(\Omega_i)-H^2(\Omega_j)\\
&<\frac12H^2(\partial C)+2\d=H^2(T_C)-2\d<H^2(E\cap C),\end{split}\ee
which contradicts the fact that $E$ is minimal. As a result, if $E$ does not contain any $\T$ type point, then there is no deformation $f$ of $E$ in $C$ such that $C\bs f(E)$ contains a segment that connects two different $\Omega_i$. On the other hand, if $E$ contains a $\T$ point, the by the argument around (0.11), $E$ is in fact the $T$ centered at this $\T$ point. In this case there is no such deformation $f$, either. We have thus

\begin{pro}Let $E$ be a 2-dimensional Almgren-minimal set in $\R^3$ such that

1) $d_{0,2}(E,T)<\min\{\e_1,\e_2\}$ ;

2) $E$ does not contain any $\T$ point. 

Let $C$, $\Omega_i$ be as above, then there exists no deformation $f$ of $E$ in $C$ such that $C\bs f(E)$ contains a segment that connects two different $\Omega_i, 1\le i\le 4$.
\end{pro}

\begin{rem}
By Proposition 2.10, the tube $\cal T$ along $\gamma_1$ cannot be too simple. For example if there exists a Lipschitz homeomorphism $f$ which is a deformation in $C$ such that
\be f(\gamma_1)=[b_1,b_2],\ee
(in this case we say that $\gamma_1$ is not "knotted"), 
then
\be C\bs f(E)=f(C\bs E)\supset f(\gamma_1)=[b_1,b_2],\ee 
which contradicts Proposition 2.10.  Thus we get the following:\end{rem}

\begin{cor}If $E$ contains no $\T$ point, then both of the $\gamma_i,i=1,2$ are "knotted". \end{cor}

After this corollary, the potential counterexample $E_0$ proposed in \cite{DJT}  is not a real counterexample, since both $\gamma_i, i=1,2$ in this example are not knotted. (We'll also explain topologically how does $E_0$ look like in the next section). Thus we have

\begin{cor} The set $E_0$ given in Section 19 of \cite{DJT} is not Almgren-minimal.
\end{cor}

To sum up, if a minimal set $E$ satisfies (2.2), then both $\gamma_1$ and $\gamma_2$ are knotted. It is not easy to imagine how to knot a $\Y$ set without producing new singularities. However this kind of set does exist. We will construct an example in the next section.

\section{A set that admits two knotted $\Y$ curves}

The purpose of this section is to give a topological example of $E$, where $\gamma_1$ and $\gamma_2$ are both knotted. But let us first look at what is the well known example $E_0$, because such an example asks already for certain imagination. In this example both $\gamma_i,i=1,2$ are not knotted.

We take a torus $T_0$ (see Figure 3-1 below). Denote by $C_0$ (the green circle in the figure) the longest horizontal circle (the equator), and fix no matter which vertical circle $L_0$ in $T_0$ (the red circle in the figure). Denote by $x_0$ their intersection. Take $r_0>0$ such that $B_0=B(x_0,r_0)\cap T_0$ (the blue circle) is a non-degenerate topological disc. Denote by $a_1,a_2$ the intersection of $\partial B_0$ and $C_0$, and $b_1,b_2$ the intersection of $\partial B_0$ and $L_0$.

\centerline{\includegraphics[width=0.9\textwidth]{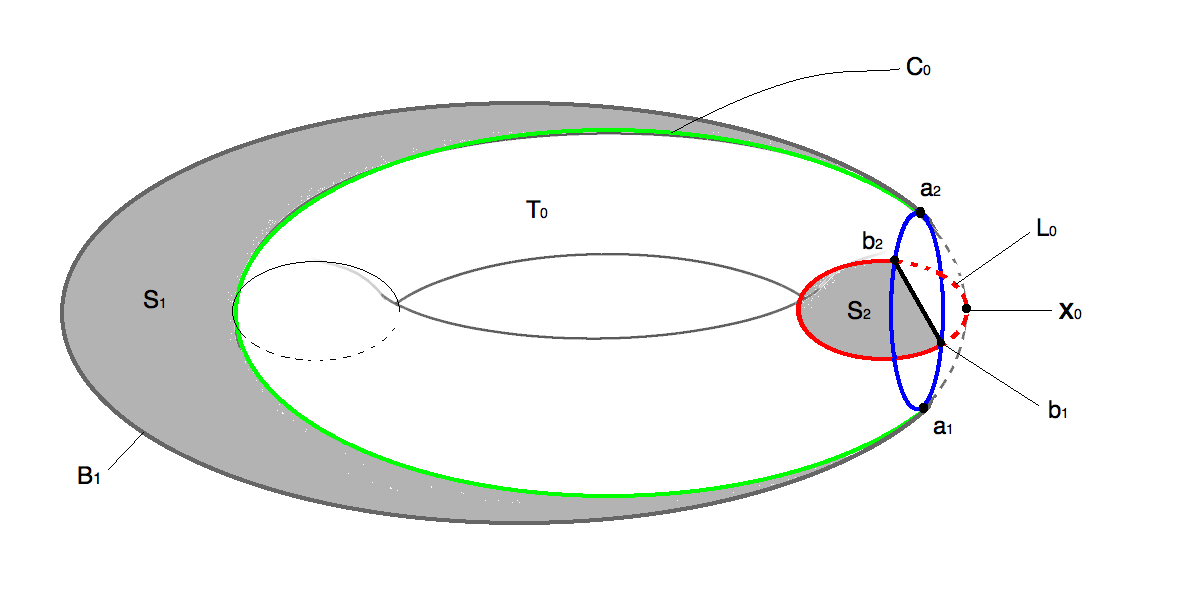}}
\nopagebreak[4]
\centerline{3-1}

Denote by $\widetilde{a_1a_2}=C_0\bs B_0$ the arc between $a_1$ and $a_2$, and $\widetilde{b_1b_2}=L_0\bs B_0$ the arc between $b_1$ and $b_2$. Next denote by $S_2$ the vertical planar part whose boundary is $[b_1,b_2]\cup\widetilde{b_1b_2}$. On the other hand, denote by $P$ the plane containing $C_0$, and take a closed disc $\overline B_1\subset P$ which contains $\widetilde{a_1a_2}$ and whose boundary contains $a_1$ and $a_2$. Now denote by $\widehat{a_1a_2}=\partial B_1\bs B_0$ the larger arc of $\partial B_1$ between $a_1$ and $a_2$, and denote by $S_1\subset P$ the part between $\widetilde{a_1a_2}$ and $\widehat{a_1a_2}$.

Now we are happy to announce that the set $(T_0\bs B_0)\cup S_1\cup S_2$ is topologically the example $E_0$ given in Section 19 of \cite{DJT}. Here $a_1,a_2,b_1,b_2$ are the four $\Y$ points which correspond to the four $\Y$ points in $E_0\cap\partial B(0,1)$, $\widetilde{a_1a_2}$ and $\widetilde{b_1b_2}$ correspond to $\gamma_1$ and $\gamma_2$ respectively; $\widehat{a_1a_2}$ and $[b_1,b_2]$, together with the four arcs on $\partial B_0$ between $a_i$ and $b_j,i,j=1,2$ correspond to the six curves of $E_0\cap\partial B(0,1)$. And if we tried to modify our topological example so that the surfaces meet each other with $120^\circ$ angles along the curves $\gamma_i$, then we would get $E_0$.


\smallskip

After the above discussion, we are now ready to construct (in $\R^3$) our example $E_1$ where $\gamma_1$ and $\gamma_2$ are both knotted. Moreover, in $\R^3\bs E_1$ there is no non-knotted curve that connects $\Omega_1$ to $\Omega_2$, or $\Omega_3$ to $\Omega_4$. The idea is to replace the $\gamma_1$, $\gamma_2$ in $E_0$, which are a pair of co-generators of $\pi_1(T_0)$, by another pair of knotted represents of co-generators of $\pi_1(T_0)$ torus $T_0$.

Let us first point out that the following example $E_1$ is just a topological one, and it is not very likely that $E_1$ could be minimal.

Still take our torus $T_0$. Take all the notation as before, that is, $C_0$ denotes the longest horizontal circle (the equator), $L_0$ is a vertical circle in $T_0$. Denote by $x_0$ their intersection. Take $r_0>0$ such that $B_0=B(x_0,r_0)\cap T_0$ (the blue circle) is a non-degenerate topological disc. Denote by $a_1,a_2$ the intersection of $\partial B_0$ and $C_0$, and $b_1,b_2$ the intersection of $\partial B_0$ and $L_0$.

Denote by $\Gamma=\Z^2$ the integer lattice in $\R^2$. We identify $T_0$ with the image of $\pi:\R^2\to \R^2/\Z^2$. For any two integers $m,n\in \Z$, denote by $d(m,n)=[(0,0),(m,n)]$ the segments of endpoints $(0,0)$ and $(m,n)$. Denote by $k(m,n)=\pi(d(m,n))$. Then $k(m,n)$ is a simple closed curve (that is, $\pi$ is injective on $d(m,n)$) if and only if their greatest common divisor $(m,n)=1$. For any integers $m,n,a,b$ with $(m,n)=(a,b)=1$, $K(m,n)$ and $K(a,b)$ represent a pair of co-generators of $\pi_1(T_0)$ if and only if $|det\left(\begin{array}{cc}m&n\\a&b\end{array}\right)|=1$. (c.f.\cite{Rol}). Without loss of generality, suppose that $K(1,0)=L_0$, and $K(0,1)=C_0$.

Take $K(2,3)$ and $K(3,4)$ a pair of knotted curves which represent  a pair of co-generators of $\pi_1(\T_0)$. Then the two curves intersect each other at one point. Without loss of generality, suppose this point of intersection is $x_0$. Denote by $Int(T_0)$ and $Ext(T_0)$ the two connected components of $S^3\bs T_0$.

First we want to construct two surfaces $S_1$ and $S_2$, such that $S_1\subset Ext(T_0), S_2\subset Int(T_0)$, and $\partial S_1=K(2,3), \partial S_2=K(3,4)$.

Notice that the torus knot $K(3,2)$ is a trefoil knot (see Figure 3-2 left), which bounds a non orientable surface $S_1'\subset Int (T_0)$ (See Figure 3-2 right). The pair of topological spaces $(Int(T_0)\cup T_0,T_0)$ is homeomorphic to $(Ext(T_0)\cup T_0, T_0)$, by some homeomorphism $\varphi$ or $S^3$ that sends the point $\infty$ to a point in $Int(T_0)$, and $\varphi(K(3,2))=K(2,3)$. Thus $K(2,3)$ bounds a surface $S_1=\varphi(S_1')\subset Ext(T_0)$.

\centerline{\includegraphics[width=0.3\textwidth]{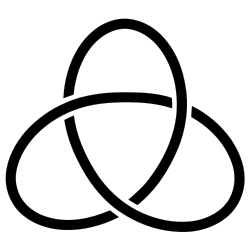}\hskip 1cm\includegraphics[width=0.32\textwidth]{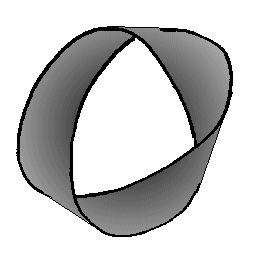}}
\nopagebreak[4]
\centerline{3-2}

For the curve $K(3,4)$, it intersects with the vertical circle $L_0$ at four points $p_0,p_1,p_2,p_3$ in clockwise order. Denote by $s_0\subset Int(T_0)$ the vertical planar disk whose boundary is $L_0$. And for $\theta\in[0,2\pi]$, $s_\theta$ denote the vertical section disk of $Int(T_0)$ with polar angle $\theta$. (See Figure 3-3 below). Then $s_\theta$ also intersects $K(3,4)$ at four points, and when $\theta<2\pi$, the intersection of $K(3,4)$ with the tube $\cup_{0\le \a\le \theta}s_\a$ is the disjoint union of four curves. Then for each $0\le i\le 3$, there is a point on $K(3,4)\cap s_\theta$ that is connected to $p_i$ by one of these four curves. Denote by this point $p_i(\theta)$. Notice that at the angle $2\pi$, we have $s_0=s_{2\pi}$, hence for each $0\le i\le 3$, $p_i(2\pi)$ is one of the two points on $s_0$ that is adjacent to $p_i$, and $\{p_i(2\pi),0\le i\le 3\}=\{p_i(0),0\le i\le 3\}$. Moreover, since $K(3,4)=\cup_{0\le i\le 3}p_i([0,2\pi[)$ is connected, if $p_i(2\pi)=p_j$, then $p_j(2\pi)\ne p_i$. Under all these conditions, there are only two possibilities : either $p_i(2\pi)=p_{i+1}$ for $i=0,1,2$ and $p_3(2\pi)=p_0$, either $p_i(2\pi)=p_{i-1}$ for $i=1,2,3$ and $p_0(2\pi)=p_3$. Without loss of generality, suppose that $p_i(2\pi)=p_{i+1}$ for $i=0,1,2$ and $p_3(2\pi)=p_0$.

\centerline{\includegraphics[width=0.8\textwidth]{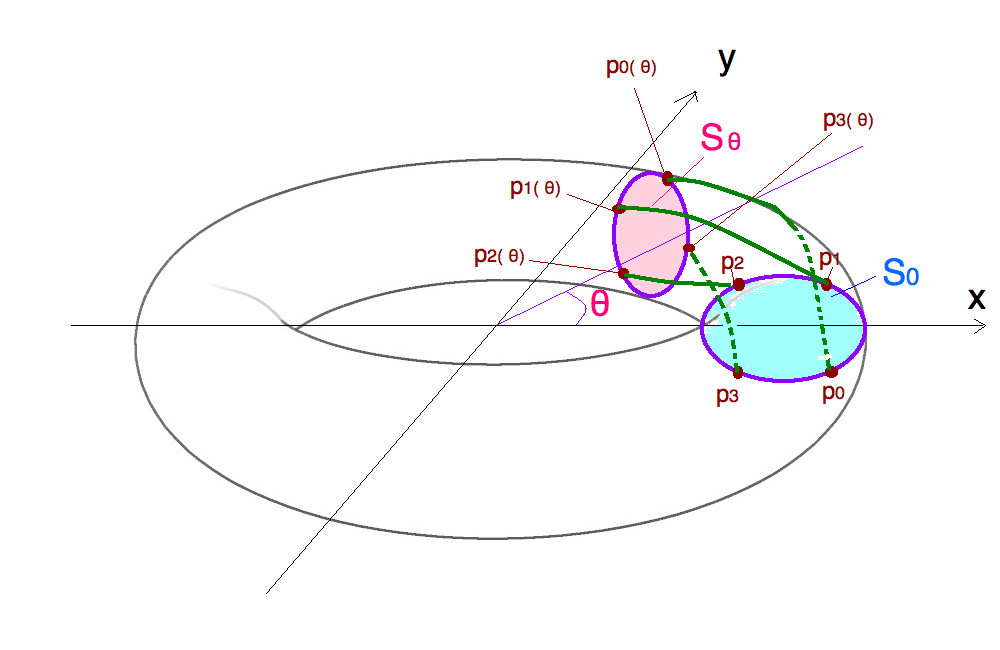}}
\nopagebreak[4]
\centerline{3-3}

Take an isotopy $f:[0,2\pi]\times s_0\to Int(T_0)\cup T_0$, such that $f_0=id$, $f_\theta(s_0)=f(\theta, s_0)=s_\theta$, $f_\theta(\partial s_0)=\partial s_\theta$, and $f_\theta(p_i)=p_i(\theta)$. Then the image $S_2'=f([0,\frac32\pi]\times ([p_0,p_1]\cup[p_2,p_3]))$ is a surface inside $T_0$, whose boundary is the curve $\{K(3,4)\cap (\cup_{0\le\theta\le\frac32\pi}s_\theta)\}\cup([p_0,p_1]\cup[p_2,p_3])\cup f_{\frac32\pi}([p_0,p_1]\cup[p_2,p_3])$.

Now we have to find a surface in the rest part, i.e. in $\cup_{\frac32\pi\le\theta\le 2\pi}s_\theta$, whose boundary is $\{K(3,4)\cap (\cup_{\frac32\pi\le\theta\le2\pi}s_\theta)\}\cup([p_0,p_1]\cup[p_2,p_3])\cup f_{\frac32\pi}([p_0,p_1]\cup[p_2,p_3])$. Notice that we cannot continue to use the image by $f_t, \frac32\pi\le\theta\le2\pi$, because $f_{2\pi}([p_0,p_1])$ will be something that connects $p_1$ and $p_2$, rather than a curve that connects $p_0$ to $p_1$ or $p_2$ to $p_3$. But we can find the solution by a saddle surface $S_2''$. Refer to Figure 3-4 below, where $a_i$ denotes $f_{\frac32\pi}(p_i)$.

\centerline{\includegraphics[width=0.5\textwidth]{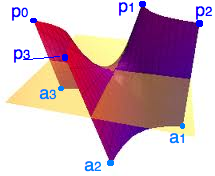}}
\nopagebreak[4]
\centerline{3-4}

Denote by $S_2=S_2'\cup S_2''\subset Int(T_0)$, then $\partial S_2=K(3,4).$

Now to sum up, we have found two surfaces $S_1\in Ext(T_0)$, $S_2\in Int(T_0)$, with $\partial S_1=K(2,3)$, $\partial S_2=K(3,4)$.

Now we take a diffeomorphism of $S^3$, which maps $T_0$ to $T_0$, $Int(T_0)$ to $Int(T_0)$, and $Ext(T_0)$ to $Ext(T_0)$. Moreover we ask that the images $l_1$ of $K(2,3)$ and $l_2$ of $K(3,4)$ satisfies that $l_1\cap l_2=x_0$, $l_1\cap B_0=C_0\cap B_0$ the shorter arc of $C_0$ between $a_1$ and $a_2$, and $l_2 \cap B_0=L_0\cap B_0$ the arc of $L_0$ between $b_1$ and $b_2$ that passes through $x_0$. Then the images of $S_1$ and $S_2$ are still two surfaces $S_3$, $S_4$, with $\partial S_3=l_1,\partial S_4=l_2$. $S_3\subset Ext(T_0), S_4\subset Int(T_0)$.

We still need a little modification, because, the two surfaces $S_3$ and $S_4$ meet each other at the boundary. So we take a homeomorphism $\varphi$ of $S^3$, which fix $T_0\bs B_0$, and $\varphi (l_1\cap B_0)=\widehat{a_1,a_2}$, $\varphi(l_2\cap B_0)=[b_1,b_2]$. $S_5=\varphi(S_3)\subset Ext(T_0)$, and $S_6=\varphi(S_4)$ is contained in $Int(T_0)\bs C$, where $C$ denotes the convex hull of $\{a_1,a_2,b_1,b_2,x_0\}$. 

Denote by $E_1=(T_0\bs B_0)\cup S_5\cup S_6$. Denote by $\gamma_1=l_1\bs B_0$, and $\gamma_2=l_2\bs B_0$. These are two knotted $\Y$ curves of $E_1$, because as in $E_0$, we have two surfaces $S_5,S_6$, such that $\partial S_5=\gamma_1\cup \widehat{a_1,a_2}$ and $\partial S_6=\gamma_2\cup[b_1,b_2]$. We can deform $E_1$ into $B(0,1)$, such that $E_1\cap\partial B(0,1)=T\cap\partial B(0,1)$ for some $\T$ set $T$. Here $a_1,a_2,b_1,b_2$ are the four $\Y$ points which correspond to the four $\Y$ points in $E_1\cap\partial B(0,1)$, $\widehat{a_1a_2}$ and $[b_1,b_2]$, together with the four arcs on $\partial B_0$ between $a_i$ and $b_j,i,j=1,2$ correspond to the six curves of $E_1\cap\partial B(0,1)$. 
Thus we have constructed an example whose set of $\Y$ points is the union of two knotted curves, moreover we cannot find any non-knotted curves in $B(0,1)\bs E_1$ that connects $a_1$ to $a_2$, or $b_1$ to $b_2$. That is, there is no deformation $f$ of $E_1$ in $B(0,1)$ such that $B(0,1)\bs f(E_1)$ contains a segment that connects two different $\Omega_i$. But this example $E_1$ seems too complicated to be minimal. We do not know how to prove this. 

However for another closely related problem, we can prove that a minimal set in that case could not admit a knotted $\Y$ curve. See the next section.

\section{Another related problem}

We take a $\T$ set $T$ centered at the origin. That is, $T$ is the cone over the 1-skeleton of a regular tetrahedron $C$ centered at the origin and inscribed to the unit ball.

In this section we will discuss whether there exists a set $E\subset\overline B(0,1)$ different from $T\cap\overline B(0,1)$, which is minimal in $B(0,1)$, and such that $E\cap\partial B(0,1)=T\cap\partial B(0,1)$.

Denote by $B=B(0,1)\subset\R^3$, and $\overline B$ its closure. Then $T$ divides the sphere $\partial B$ into four equal triangular open regions $\{S_i\}_{1\le i\le 4}$, with 
\be \cup_{i=1}^4 \overline{S_i}=\partial B\mbox{ and }\cup_{i=1}^4 S_i=\partial B\bs T.\ee
Recall that $T$ divides $\partial C$ into four equal open planar triangles $\{\Omega_i\}_{1\le i\le 4}$. For notation convenience we ask that for each $i$, $S_i$ and $\Omega_i$ share the same three vertices.

Denote by $a_j,1\le j\le 4$ the four vertices of $T\cap\partial B$, where $a_j=\cap_{i\ne j}\overline{S_j}\cap\partial B$ is the point opposite to $S_j$.

\begin{pro}

Let $E\subset \overline B\cap \R^3$ be a closed, 2-rectifiable, locally Ahlfors regular set, with
\be E\cap \partial B=T\cap\partial B.\ee

Then 

1) If $H^2(E)<H^2(T\cap B)$, 
\be \begin{split}&\mbox{there exists }1\le i<j\le 4\mbox{, and four points }a,b,c,d\mbox{ which lie in a common plane},\\
&\mbox{ such that }a\in S_i, d\in S_j,b,c\in B \bs E,\angle abc>\frac\pi 2,\angle bcd>\frac\pi2\\
&\mbox{ and }[a,b]\cup [b,c]\cup [c,d]\subset \overline B\bs E.\end{split}\ee
Here $[x,y]$ denotes the segment of endpoints $x$ and $y$, and $\angle abc\in[0,\pi]$ denotes the angle of the smaller sector bounded by $\overline{ba}$ and $\overline {bc}$.

2)If $E$ is a reduced minimal set in $B$ and verifies (4.3), then
\be \mbox{ either }E=T\cap \overline B\mbox{, or (4.4) is true.}\ee

\end{pro}

Before we prove Proposition 4.2, let us first give a corollary.

\begin{cor}Let $E\subset\overline B$ be a reduced minimal set in $B$ and satisfy (4.3). Then if $E\ne T\cap\overline B$, we have
\be H^2(E)\le H^2((T\cap\overline B\bs C)\cup G)=H^2(T\cap\overline B)-(2\sqrt 2-\frac43\sqrt 3)(\approx 0.519).\ee
\end{cor}

\nd Let $E$ be such a set. Then by (4.5), (4.4) is true. But (4.4) gives the existence of a deformation $f$ in $B$ such that $f(E)\subset B\bs [a,d]$, hence we can deform $E$ on a subset of $(T\cap B\bs C^\circ)\cup G$ ($C^\circ$ denotes the interior of $C$), where $G=\partial C\bs(\Omega_i\cup\Omega_j)$ (recall that $\Omega_i,\Omega_j$ are the two faces of $C$ corresponding to the two faces $S_i,S_j$ of $\partial B$, where $S_i,S_j$ contain the points $a$ and $d$).

Thus by (2.4), 
\be H^2(E)\le H^2((T\cap \overline B\bs C)\cup G)=H^2(T\cap\overline B)-(2\sqrt 2-\frac43\sqrt 3).\ee\qed

\noindent Proof of Proposition 4.2.

We are going to prove 1) by contraposition. Suppose that (4.4) is not true.

Denote by $P_j$ the plane orthogonal to $\overrightarrow{oa_j}$ and tangent to the unit sphere, and denote by $p_j$ the orthogonal projection to $P_j$. Set $R_j=p_j(\overline S_j)\subset p_j(\cup_{i\ne j}\overline S_i)\subset P_j$. Then for each $1\le j\le 4$ and each $x\in R_j$,
\be p_j^{-1}(x)\cap E\ne\emptyset.\ee

In fact if (4.9) is not true for some $j$, that is, $R_j\bs p_j(E)\ne\emptyset$. As the projection of a compact set, $p_j(E)$ is compact. Thus $R_j\bs p_j(E)$ is a non empty open set. Note that $R_j\bs(\cup_{i\ne j}p_j(S_i))$ is of measure zero, therefore
\be (R_j\bs p_j(E))\cap (\cup_{i\ne j}p_j(S_i))\ne\emptyset.\ee

Take $x\in (R_j\bs p_j(E))\cap (\cup_{i\ne j}p_j(S_i))$. Then $x\not\in\partial R_j$, because $\partial S_j\subset E$ and hence $\partial R_j=\partial p_j(S_j)=p_j(\partial S_j)\subset p_j(E)$. As a result, $p_j^{-1}(x)\cap B$ is a segment $[a,d]$ perpendicular to $P_j$, with $a\ne d$, $a\in S_j^\circ$ and $d\in\cup_{i\ne j} S_i^\circ$.  Take $b,c\in[a,d]$ such that $a,b,c,d$ are different. Then (4.4) holds, which contradicts our hypothesis.

Hence (4.9) holds. Now for each $x\in R_j$, denote by $f_j(x)$ the point in $p_j^{-1}(x)\cap E$ which is the nearest one to $R_j$. In other words, $f_j(x)$ is the first point in $E$ whose projection is $x$. This point exists by (4.9), and is unique, since $p_j^{-1}(x)$ is a line orthogonal to $R_j$.

Denote by $A_j=f_j(R_j)$. Then $A_j$ is measurable. In fact, 
\be \begin{split}A_j&=\{x\in E:\forall y\in E\mbox{ such that }d(y,P_j)<d(x, P_j),|p_j(y)-p_j(x)|>0\}\\
&=\bigcap_{p,q}\{x\in E:\forall y\in E\mbox{ such that }d(y,P_j)<d(x, P_j)-2^{-p},|p_j(y)-p_j(x)|>2^{-q}\}
.\end{split}\ee

Now $E$ is rectifiable, hence $A_j\subset E$ is also rectifiable. Therefore for almost all $x\in A_j$, the approximate tangent plane $T_xA_j$ of $A_j$ at $x$ exists. Denote by $v_j=\frac{\overrightarrow{oa_j}}{|oa_j|}$ the unit exterior normal vector of $P_j$, and denote by $w_j(x)$ the unit vector orthogonal to $T_xA_j$ such that $<v_j,w_j(x)>\ge 0$. Then $w_j(x)$ is well defined for every $x\in A_j$ with $T_xA_j\not\perp P_j$.

Denote by 
\be E_j=\{x\in A_j:T_xA_j\not\perp P_j\}.\ee
Then $w_j$ is a measurable vector field on $E_j$.
On the other hand, by Sard's theorem, $H^2(p_j(A_j\bs E_j))=0.$
But $p_j$ is injective on $A_j$, hence $p_j(A_j\bs E_j)=p_j(A_j)\bs p_j(E_j)=R_j\bs p_j(E_j)$, and thus 
\be H^2(R_j\bs p_j(E_j))=0.\ee Moreover, for almost all $x\in E_j, T_xA_j=T_xE_j$.

We are going to show that
\be \int_{E_j} <v_j,w_j(x)> dx=H^2(R_j).\ee

First, we apply the area formula for Lipschitz maps between rectifiable sets in \cite{Fe} 3.2.20, with $m=\nu=2$, $W=E_j$, $f=p_j,g=1_{R_j}$, and we get
\be \int_{E_j} ||\wedge_2 ap Dp_j(x)||dH^2 x=\int_{R_j}N(p_j,z)dH^2z.\ee

Moreover by (4.9) and (4.13), $N(p_j,z)\ge 1$ for almost all $z\in R_j$. On the other hand, $N(p_j,z)\le 1$ since $E_j$ is contained in the set $A_j$ on which $p_j$ is injective. Hence $N(p_j,z)=1$ for almost all $z\in R_j$. Therefore
\be \int_{R_j}N(p_j,z)dH^2z=H^2(R_j).\ee
For the left side of (4.15), take $w_j^1(x)$ a unit vector in $T_xE_j$ such that $w_j^1(x) // R_j$, and $w_j^2(x)$ the unit vector in $T_xE_j$ which is orthogonal to $w_j^1(x)$. Then $p_j(w_j^1(x))\perp p_j(w_j^2(x))$, by elementary geometry in $\R^3$. Therefore
\be \begin{split}||\wedge_2 ap Dp_j(x)||&=||p_j(w_j^1(x))\wedge p_j^j(w_2(x))||=|p_j(w_j^1(x))||p_j(w_j^2(x))|\\
&=|p_j(w_j^2(x))|.\end{split}\ee
The first inequality is because $p_j$ is a linear map from $\R^2$ to $\R^2$, the second inequality is because $p_j(w_j^1(x))\perp p_j(w_j^2(x))$, and the last is because $w_j^1(x) // S_j$.

Now set $v_j^2(x)=\frac{p_j(w_j^2(x))}{|p_j(w_j^2(x))|}\in P_j$. This is well defined because $T_xE_j\not\perp P_j$ and hence $|p_j(w_j^2(x))|>0$. Then $w_j(x), w_j^2(x), v_j, v_j^2(x)$ are all orthogonal to $w_j^1(x)$, and hence belong to a same plane, with $w_j(x)\perp w_j^2(x), v_j\perp v_j^2(x)$. Therefore 
\be |<w_j(x),v_j>|=|<w_j^2(x),v_j^2(x)>|=|p_j(w_j^2(x))|.\ee
But by definition, $<w_j(x),v_j>\ge 0$, hence
\be <w_j(x),v_j>=|p_j(w_j^2(x))|=||\wedge_2 ap Dp_j(x)||\ee
by (4.17). Combine (4.15), (4.16) and (4.19), we get (4.14). Note that $v_j$ does not depend on $E$. 

Now for $x\in A_j\bs E_j$, we define a measurable vector field $w_j(x)$ such that $w_j(x)\perp T_xA_j$. Then $<w_j(x),v_j>=0$ for almost all $x\in A_j\bs E_j$. Hence we have
\be \int_{A_j} <v_j,w_j(x)> dx=H^2(R_j).\ee

We sum over $j$, and get 
\be \sum_{j=1}^4 \int_{A_j} <v_j,w_j(x)> dx=\sum_{j=1}^4 \int_{E_j} <v_j,w_j(x)> dx=\sum_{i=1}^4 H^2(R_j).\ee

Next, set $E_j^0=E_j\bs\cup_{i\ne j}E_i$, $E_{ij}=(E_i\cap E_j)\bs\cup_{k\ne i,j}E_k$ for $i\ne j$. We claim that
\be E_j\bs (E_j^0\cup \cup_{i\ne j}E_{ij})\mbox{ is of measure zero for all }j.\ee

Suppose (4.22) is not true. Then there exists three different $i,j,k$ such that $E_i\cap E_j\cap E_k$ is of positive measure. Suppose for example that $i=1,j=2,k=3$, and set $E_{123}=E_1\cap E_2\cap E_3$. Now since $E_{123}$ is a measurable rectifiable set of positive measure,  and that $E_{123}\subset E$, for almost all $x\in E_{123}$, the approximate tangent plane $T_xE_{123}$ of $E_{123}$ at $x$ exists and is equal to $T_xE$. Moreover since $E$ is locally Ahlfors regular, $T_xE$ is a real tangent plane (c.f. for example \cite{DMS}, Exercise 41.21, page 277). We choose and fix a such $x\in E_{123}$.

By definition of $A_j$, for $j=1,2,3$, the segment $[x,p_j(x)]\cap E=\{x\}$. And by definition of $E_j$, $T_xE\not\perp P_j$, and hence $[x,p_j(x)]\cap (T_xE+x)=\{x\}$, since $[x,p_j(x)]\perp P_j$. The affine subspace $T_xE+x$ separates $\R^3$ into 2 half spaces, and since for $j=1,2,3,]x,p_j(x)]\cap (T_xE+x)=\emptyset$, there exist $1\le i<j\le 3$ such that $]x,p_i(x)]$ and $]x,p_j(x)]$ are on the same side of $T_xE+x$. Suppose for example that $i=1,j=2$.

For $i=1,2$, denote by $\a_i$ the angle between $[x,p_i(x)]$ and $T_xE+x$. Set $\a=\min\{\a_1,\a_2\}$. Then since $T_xE$ is a real tangent plane, there exists $r>0$ such that for all $y\in E\cap B(x,r)$,
\be d(y,T_xE+x)< \frac r2\sin\a .\ee 

Set $b=[x,p_1(x)]\cap \partial B(x,r), c=[x,p_2(x)]\cap \partial B(x,r)$. Then by definition of $\a$, $d(b,T_xE+x)\ge r\sin\a, $ and $d(c,T_xE+x)\ge r\sin\a$. But $b,c$ are at the same side of $T_xE+x$, hence for all $y\in[b,c]$, $d(y,T_xE+x)\ge r\sin\a$, and hence $[b,c]\cap E=\emptyset$, because of (4.23). 

Now set $a=p_1(x),d=p_2(x)$. Note that in the triangle $\Delta_{xbc}$, $|xb|=|xc|$, which gives that $\angle xbc=\angle xcb$. But $\angle xbc+\angle xcb+\angle bxc=\pi$, $\angle bxc>0$, Hence $\angle xbc=\angle xcb<\frac\pi2$. As a result, $\angle abc=\pi-\angle xbc>\frac\pi2$, $\angle bcd=\pi -\angle xcb>\frac\pi2$. Thus we have found four points $a,b,c,d$ such that (4.4) is true, which contradicts our hypothesis.

Thus we get (4.22). And consequently we have 
\be H^2(\cup_{j=1}^4 E_j)=\sum_{j=1}^4 H^2(E_j^0)+\sum_{1\le i<j\le 4}H^2(E_{ij}).\ee

For estimating the measure, we are going to use the paired calibration method (introduced in \cite{LM94}). Recall that $v_j$ is the unit exterior normal vector of $P_j$. Thus by (4.21),
\be \begin{split}&\sum_{i=1}^4H^2(R_j)=\sum_{j=1}^4\int_{E_j}<v_j,w_j(x)>dx\\
=&\sum_{j=1}^4\int_{E_j^0}<v_j,w_j(x)>dx+\sum_{1\le i<j\le 4}\int_{E_{ij}}<v_i,w_i(x)>+<v_j,w_j(x)>dx\\
\end{split}\ee
For the first term,
\be \begin{split}|\int_{E_j^0}<v_j,w_j(x)>dx|&\le \int_{E_j^0}|<v_j,w_j(x)>|dx\\
&\le \int_{E_j^0}|v_j||w_j(x)|dx=H^2(E_j^0)\end{split}\ee
and hence
\be |\sum_{j=1}^4\int_{E_j^0}<v_j,w_j(x)>dx|\le \sum_{j=1}^4|\int_{E_j^0}<v_j,w_j(x)>dx|\le \sum_{j=1}^4 H^2(E_j^0).\ee
 For the second term, observe that $w_i(x)=\pm w_j(x)$ for $x\in E_{ij}$, hence we set $\e_x=\frac{w_i(x)}{w_j(x)}$. Then
 \be \begin{split}|<v_i,&w_i(x)>+<v_j,w_j(x)>|=|<v_i+\e(x)v_j, w_i(x)>|\\
 &\le |v_i+\e(x)v_j||w_i(x)|=|v_i+\e(x)v_j|\le\max\{|v_i+v_j|,|v_i-v_j|\}.\end{split}\ee
 
 By definition of $v_j$, the angle between $v_i$ and $v_j$ is the supplementary angle of the angle $\theta_{ij}$ between $P_i$ et $P_j$. Thus a simple calculus gives
 \be |v_i+v_j|=\frac{2}{\sqrt 3}<1, |v_i-v_j|=\frac{2\sqrt 2}{\sqrt 3}>1.\ee
Hence $ \max\{|v_i+v_j|,|v_i-v_j|\}=|v_i-v_j|>1$. Denote by  $D$ this value. By (4.28), 
\be \begin{split}|\int_{E_{ij}}<v_i,&w_i(x)>+<v_j,w_j(x)>dx|\\
&\le\int_{E_{ij}}|<v_i,w_i(x)>+<v_j,w_j(x)>|dx\le DH^2(E_{ij})\end{split}\ee
and hence
\be \begin{split}|&\sum_{1\le i<j\le 4}\int_{E_{ij}}<v_i,w_i(x)>+<v_j,w_j(x)>dx|\\
&\le \sum_{1\le i<j\le 4}|\int_{E_{ij}}<v_i,w_i(x)>+<v_j,w_j(x)>dx|=D\sum_{1\le i<j\le 4}H^2(E_{ij}).\end{split}\ee

Combine (4.25), (4.27) and (4.31), we get
\be \begin{split}\sum_{i=1}^4H^2(R_j)&\le \sum_{j=1}^4 H^2(E_j^0)+D\sum_{1\le i<j\le 4}H^2(E_{ij})\\
&\le D[\sum_{j=1}^4 H^2(E_j^0)+\sum_{1\le i<j\le 4}H^2(E_{ij})]\mbox{ (since }D>1)\\
&=D H^2(\cup_{j=1}^4 E_j)\le DH^2(E).\end{split}\ee
 
On the other hand, we can do the same thing for $T$, the cone over the 1-skeleton of the regular tetrahedron $C$. Since $T$ separates the four faces of $C$, (4.4) is automatically false for $T$. Then by the notions above, we can see that $T_i^0=\emptyset$ for all $i$, $\e_{ij}=-1$ for all $i\ne j$, and $(v_i-v_j)\perp T_xT$ for almost all $x\in T_{ij}$, which implies that
\be <v_i,w_i(x)>+<v_j,w_j(x)>=D\ee
for all $x\in T_{ij}$. So briefly, all the inequalities in the whole argument above are equalities for  $T\cap\overline B$. As a result,
\be DH^2(E)\ge \sum_{i=1}^4H^2(R_j)= DH^2(T\cap\overline B),\ee
and hence
\be H^2(E)\ge H^2(T\cap\overline B)\ee for all $E$ that does not verify (4.4).

\bigskip

Now let us prove 2). Let $E$ be a reduced minimal set, then it is rectifiable and locally Ahlfors regular in $B$ (c.f.\cite{DS00}). 

First note that $H^2(E)\le H^2(T\cap\overline B)$. In fact, for each $x\in \overline B\bs T$, there exists $1\le i\le 4$, such that $x$ and $S_i$ belong to the same connected component of $B\bs T$. Then denote by $f(x)$ the first intersection of $x+[0,a_i)$ with $T$. Then $f:\overline B\to T$ is a 2-Lipschitz retraction (see Figure 4-1). Now if $E$ is a minimal set that verifies (4.3), for each $\e>0$, we define $g_\e:\partial B\cup E\to T\cup \partial B$ by $g(x)=f(x)$ for $x\in E\cap \overline B(0,1-\e)$, $g(x)=x$ for $x\in\partial B$. Then we can extend $g_\e$ to a 2-Lipschitz map which sends $\overline B$ to $\overline B$, by Kirszbraun's Theorem (c.f.\cite{Fe}, Thm 2.10.43). (Figure 4-1). Thus $g_\e$ deforms $E\cap \overline B(0,1-\e)$ to a subset of $T\cap\overline B(0,1-\e)$. Thus we have
\be \begin{split}H^2(g_\e(E))&=H^2(g_\e(E\cap\overline B(0,1-\e)))+H^2(g_\e(E\bs\overline B(0,1-\e)))\\
&\le H^2(g_\e(T\cap\overline B(0,1-\e)))+Lip(g_\e)^2H^2(E\bs\overline B(0,1-\e)\\
&=H^2(T\cap\overline B(0,1-\e))+4H^2(E\bs\overline B(0,1-\e))\\
&<H^2(T\cap\overline B)+4H^2(E\bs\overline B(0,1-\e)
\end{split}\ee
The second term tends to 0 when $\e$ tends to $0$. That is, for any $\d>0$, there exists $\e(\d)>0$ such that
\be H^2(g_\e(E))<H^2(T\cap\overline B)+\d.\ee

\centerline{\includegraphics[width=0.5\textwidth]{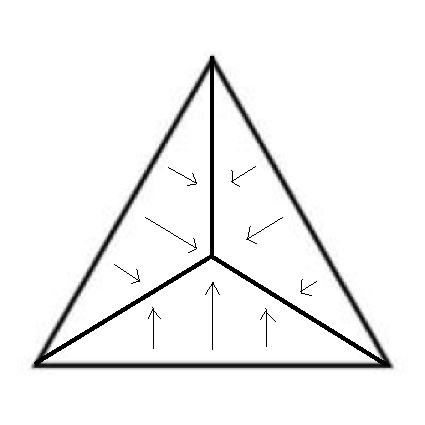}}
\nopagebreak[4]
\centerline{4-1}

Now $E$ is minimal, hence for any $\d>0$,
\be H^2(E)\le H^2(g_\e(\d)(E))\le H^2(T\cap\overline B)+\d,\ee
therefore
\be H^2(E)\le H^2(T\cap\overline B).\ee

Hence to prove 2), it is enough to prove that if (4.4) does not hold, and $H^2(E)=H^2(T\cap \overline B)$, then $E=T\cap\overline B$. In particular $E$ contains a point of type $\T$. 

By the arguments in 1), if (4.4) is not true, and $H^2(E)=H^2(T\cap\overline B)$, then the inequalities (4.26)-(4.28) et (4.30)-(4.32) are all equalities.
Thus we have

1) For almost all $x\in E_{ij}$, $T_xE_{ij}\perp v_i-v_j$. Denote by $P_{ij}$ the plane perpendicular to $v_i-v_j$, then for almost all $x\in E_{ij}$, $T_xE_{ij}=P_{ij}$.

2) For all $j$, $H^2(E_j^0)=0$, since $D>1$.

3) For all $j$, $H^2(A_j\bs E_j)=0$.

4) For all $j$, $p_j(E)=p_j(E_j)=R_j$.

Thus for almost all $x\in E$, $T_xE$ exists and is one of the $P_{ij}$. If $x$ is a point such that $T_xE$ exists, by the $C^1$ regularity (c.f.\cite {DEpi}, Theorem 1.15 and Lemma 14.4), there exists $r=r(x)>0$ such that in $B(x,r)$, $E$ is the graph of a $C^1$ function from $T_xE$ to $T_xE^\perp$, which implies that in $B(x,r)$, the function
$f: E\cap B(x,r)\to G(3,2)$, $f(y)=T_yE$ is continuous. But for $T_yE$ we have only six choices $P_{ij},1\le i<j\le 4$, which are isolated points in $G(3,2)$, and hence $T_yE=T_xE$ for all $y\in B(x,r)\cap E$. As a result $E\cap B(x,r)=(T_xE+x)\cap B(x,r)$, a disc parallel to $P_{ij}$.

Still by the $C^1$ regularity, the set $E_P=\{x\in E\cap B:T_xE\mbox{ exists}\}$ is a $C^1$ manifold, and is open in $E$. Thus we deduce that 
\be\begin{split} \mbox{each connected component of }E_P\mbox{ is part of}\\ \mbox{ a plane that is parallel to one of the }P_{ij}.\end{split}\ee

Set $E_Y=\{x\in E:x\mbox{ is of type }\Y\}$. Then $E_Y\ne\emptyset$, because otherwise by (4.40), $E\cap B$ is the intersection of $B$ with a translation of one of the $P_{ij}$, but then $E\cap\partial B$ is surely not $T\cap\partial B$.

Now if $x\in E_Y$, by the $C^1$ regularity around points of type $\Y$ (c.f.\cite {DEpi}, Theorem 1.15 and Lemma 14.6), there exists $r=r(x)>0$ such that in $B(x,r)$, $E$ is $C^1$ equivalent to a $\Y$ set $Y$. Denote by $L_Y$ the spine of $Y$, and $S_1,S_2,S_3$ the three open half planes of $Y$. Then if we denote by $\varphi$ the $C^1$ diffeomorphism which sends $Y$ onto $E$ in $B(x,r)$, the $\varphi(S_i)\cap B(x,r),1\le i\le 3$ are connected $C^1$ manifolds, and hence each of them is a part of a plane parallel to $P_{ij}$. Consequently, $\varphi(L_y)\cap B(x,r)$ is an open segment passing through $x$ and parallel to one of the $D_j,1\le j\le 4$, where $D_j=P_{ij}\cap P_{jk}$. 

Hence $E_Y\cap B$ is a union of open segments $I_1,I_2\cdots$, each of them is parallel to one of the $D_j$, and every endpoint is either a point in the sphere $\partial B$, or a point of type $T$. Moreover
\be \begin{split}\mbox{for each }x\in E_Y\mbox{ such that }T_xE_Y=D_j\mbox{, there exists }r>0\\
\mbox{ such that in }B(x,r),E\mbox{ is a  }\Y\mbox{ set whose spine is }x+D_j.
    \end{split}
\ee

Now if $x\in E$ is a $\T$ point, then by the arguments above, the blow-up limit $C_xE$ of $E$ at $x$ is the set $T$ (the set $T$ that defined at the very beginning of this section). As a result, for each segment $I_i$, at least one of its endpoints is in the unit sphere. In fact, if both of the two endpoints $x,y$ of $I_i$ are of type $\T$, then at least one of the two  blow-up limit $C_xE$ and $C_yE$ are not the set $T$, because two parallel $\T$ sets cannot be connected by a common spine.

Hence all the segments $I_i$ touch the boundary. 

\begin{lem} If $x$ is a $\T$ point (and hence $C_xE=T$), then $(T+x)\cap B\subset E$.
\end{lem}
 
\nd  By the $C^1$ regularity around points of type $T$, there exists $r>0$ such that in $B(x,r)$, $E$ is a $C^1$ version of $T+x$. Then by (4.40) and (4.41), $E\cap B(x,r)=(T+x)\cap B(x,r)$. Denote by $L_i,1\le i\le 4$ the four spines of $T+x$, then $L_i\cap B\subset E_Y$, because $L_i\cap B(x,r)$ is part of a segment $I_j\subset E_Y$, which has already an endpoint $x$ that does not belong to the unit sphere, hence the other endpoint must be in the sphere, which yields $I_j=L_i\cap B(0,1)$.

Now we take a one parameter family of open balls $B_s$ of radius $r\le s\le 1$, with $B_r=B(x,r)$ and $B_1=B(0,1)$, such that

1) $B_s\subset B_{s'}$ for all $s<s'$;

2) $\cap_{1>t>s}B_t=\overline B_s$ and $\cup_{t<s}B_t=B_s$ for all $r\le s\le 1$.

Set $R=\inf\{s>r,(T+x)\cap B_s\not\subset E\}.$ We claim that $R=1$.

Suppose this is not true. 

By definition of $B_s$, the four spines and the six faces of $T+x$ are never tangent to $\partial B_s$, $r<s<1$, since $B(x,r)\subset B_s$.  

Now for each $y\in \partial B_R\cap (T+x)$, $y$ is not a $\T$ point. In fact, if $y$ belongs to one of the $L_i$, then $y$ is a $\Y$ point, since $L_i\bs\{x\}\subset E_Y$ and $L_i\cap B\subset E_Y$; if $y$ is not a $\Y$ point, then there exists $i,j$ such that $y\in x+P_{ij}$. Thus there exists $r_y>0$ such that $B(y,r_y)\cap (x+T)$ is a disc $D_y$ centered at $y$. Now by definition of $R$, for all $s<R$, $B_s\cap (T+x)\subset E$, and hence $B_R\cap (T+x)\subset E$. Hence $D\cap B_R\cap B(y,r_y)\subset E$, which means that $y$ cannot be a point of type $\T$.

If $y$ is a point of type $\P$ (i.e., planar point), suppose for example that $y\in P_{ij}+x$. Then $T_yE=P_{ij}$. By (4.40), and since $R<1$, there exists $r_y>0$ such that $E\cap B(y,r_y)=(P_{ij}+y)\cap B(y,r_y)$. In other words,
\be \mbox{there exists }r_y>0\mbox{ such that }E\mbox{ coincides with }T+x\mbox{ in }B_R\cup B(y,r_y).\ee

If $y$ is a point of type $\Y$, then it is in one of the $L_i$. By the same argument as above, using (4.41), we get also (4.43).

Hence (4.43) is true for all $y\in \partial B_R\cap (T+x)$. But $\partial B_R\cap (T+x)$ is compact, we have thus a uniform $r>0$ such that for each $y$, (4.43) is true if we set $r_y=r$. But this contradicts the definition of $R$.

Hence $R=1$. But $B_1\subset B$ is of radius 1, hence $B_1=B$.  Then by definition of $R$ we get the conclusion of Lemma 4.42. \qed

By lemma 4.42, we know that if $x$ is a $\T$ point, then $x$ has to be the origin, because of (4.3). Hence $T\cap B\subset E$. But in this case, we have $E=T\cap\overline B$, because $H^2(E)=H^2(T\cap\overline B)$. 

We still have to discuss the case when there is no point of type $\T$. But in this case, the same kind of argument as in Lemma 4.42 gives the following.

\begin{lem}
Let $x$ be a $\Y$ point in $E$, and $T_xE_Y=D_j$. Denote by $Y_j$ the $Y$ whose spine is $D_j$. Then
\be (Y_j+x)\cap B\subset E.\ee
\end{lem}

But this is impossible, because $E\cap \partial B=T\cap \partial B$ contains no full part of $(Y_j+x)\cap\partial B(0,1)$ for any $x$ and $j$.

Hence we have $E=T\cap\overline B$. And thus (4.5). \qed

\part{Existence of a point of type $\T$ for 2-dimensional topological minimal sets in $\R^4$}

\section{Introduction}

In this part we discuss the property (0.13) for 2-dimensional topological minimal sets in $\R^4$ whose blow-in limits are $\T$ sets. This kind of sets exists trivially because a $\T$ cone is topological minimal in $\R^3$, and by Proposition 3.18 of \cite{2p}, it is topological minimal in any $\R^n$ for $n\ge 3$.

We are wondering if there is any other type of topological minimal sets in $\R^4$ that look like $\T$ sets at infinity but are not $\T$ sets themselves. Recall that in $\R^3$ there are no such sets (c.f. Proposition 18.1 of \cite{DJT}). A main useful property in $\R^3$ is that there are only two kinds of minimal cones whose density is less than that of $\T$ sets: the planes and the $\Y$ sets. Hence if the blow-in limits of a non-conical minimal set are $\T$ sets, then by the monotonicity of density, in this minimal set all points are of type $\P$ or $\Y$. And thus we can have the same properties that were stated around (2.5)-(2.7) and Figure 2-1, by the same argument.

In $\R^4$ we do not know if there exists a minimal cone whose density is between $\Y$ sets and $\T$ sets. However $\T$ sets are the only minimal cones that admit the simplest topology except for planes and $\Y$ sets. Hence it is likely that in $\R^4$ there are no minimal cones between $\Y$ sets and $\T$ sets.

So we make this additional assumption. Denote by $d_T$ the density of $\T$ sets, and we suppose that
\be \begin{split}\mbox{the only minimal cones in }&\R^4\mbox{ whose densities are less than }d_T\\
&\mbox{ are the planes and the }\Y\mbox{ sets}.\end{split}\ee 

We are going to discuss, under the assumption (5.1), the Bernstein type property for topological minimal sets in $\R^4$ that look like a $\T$ set at infinity.

\section{A topological criterion for potential counterexamples}

Throughout this section, we assume that (5.1) is true.

Let $E$ be a 2-dimensional topological minimal set in $\R^4$ that looks like a $\T$ set at infinity. That is, there exists a $\T$ set $T$ centered at the origin, and a sequence $\{r_k\}_{k\in\N}$ such that
 \be\lim_{k\to\infty}r_k\to\infty\mbox{ and }\lim_{k\to\infty}d_{0,r_k}(E,T)=0.\ee

We want to find a $\T$ type point in the set $E$.

Now the set $E$ is of codimension 2, hence the topological condition is imposed on the group $H_1(\R^4\bs E,\Z)$.

Denote by $\{y_i\}_{1\le i\le 4}$ the four $\Y$ points in $T\cap\partial B(0,1)$. Denote by $l_{ij}\subset T\cap\partial B(0,1)$ the great arc on the sphere that connects $y_i$ and $y_j$. The cone $T$ is composed of 6 closed sectors $\{T_{ij}\}_{1\le i\ne j\le 4}$, where $T_{ij}$ is the cone over $l_{ij}$. Denote by $x_{ij},1\le i\ne j\le 4$ the middle point of $l_{ij}$. Denote by $P_{ij}$ the 2-plane orthogonal to $T_{ij}$ and passing through $x_{ij}$. Set $B_{ij}=B(x_{ij},\frac {1}{10})\cap P_{ij}$, and denote by $s_{ij}$ the boundary of $B_{ij}$. Then $s_{ij}$ is a circle, that does not touch $T$, and $B_{ij}\cap T=B_{ij}\cap T_{ij}=x_{ij}$.

Fix an orthonormal basis $\{e_i\}_{1\le i\le 4}$ of $\R^4$. We are going to give an orientation to each $s_{ij}$, and denote these oriented circles by $\vec s_{ij}$.

For each $B_{ij}$, there are two orientations $\sigma_1=x\wedge y,\ \sigma_2=-x\wedge y$, where $x,y$ are two mutually orthogonal unit vectors that belong to the plane containing $B_{ij}$. Take the $k\in\{1,2\}$ such that $det_{\{e_i\}_{1\le i\le 4}}\overrightarrow{ox_{ij}}\wedge \overrightarrow{y_iy_j}\wedge \sigma_k>0$, and denote by $\overrightarrow B_{ij}$ the oriented disc $B_{ij}$ with this orientation.
Denote by $\vec s_{ij}=\partial \overrightarrow B_{ij}$ the oriented circle, and $[\vec s_{ij}]$ the element in $H_1(\R^4\bs T;\Z)$ represented by $\vec s_{ij}$. The six $[\vec s_{ij}],1\le i<j\le 4$ are all different, however they are algebraically dependent. 

Figure 6-1 gives an idea for the above definition. (But it is drawn in $\R^3$). Since $T$ is contained in $\R^3$, so if we fix an orientation for the other dimension in $\R^4$, the orientation of $B_{ij}$ defined before corresponds to one of the orientation of the line orthogonal to $T_{ij}$ in $\R^3$. And this orientation of the line corresponds to the orientation of $l_{ij}$ by the right-hand rule. Hence in Figure 6-1 we mark the orientation of $l_{ij}$ by the arrows to express the orientation of $[\vec s_{ij}]$. In the figure, the orientation $\underline{S}_{ij}$ means $[\vec s_{ij}]$.

\centerline{\includegraphics[width=0.6\textwidth]{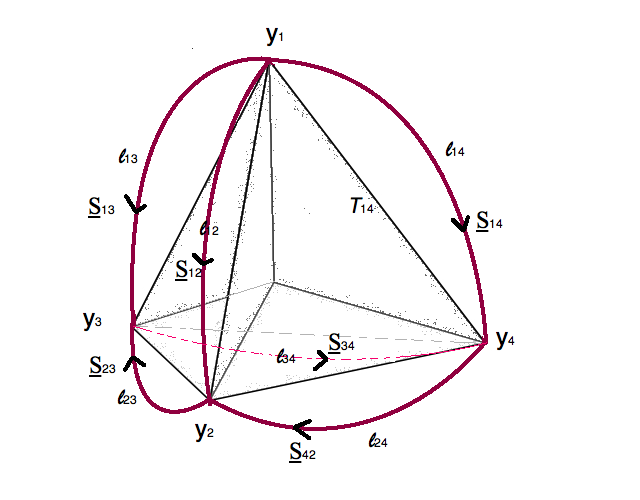}}
\nopagebreak[4]
\centerline{6-1}
\bigskip

Thus we have
\be \vec s_{ij}=-\vec s_{ji},[\vec s_{ij}]=-[\vec s_{ji}].\ee

Note that $\{[\vec s_{ij}],1\le i,j\le 4\}$ is a set of generators of the group $H_1(\R^4\backslash T;\mathbb Z)$.

We say that $[s_{ij}]$ and $[s_{kl}]$ (without the vector arrow) are different (in a homology group) if
\be [\vec s_{ij}]\ne\pm [\vec s_{kl}],\ee
and denote by $s_{ij}\sim s_{kl}$ if $[\vec s_{ij}]=\pm [\vec s_{kl}]$.

\smallskip

Return to our set $E$. Without loss of generality, we can suppose (modulo replacing $E$ by $E/r_k$ for some $k$ large) that $d_{0,3}(E,T)$ is small enough (for example less than a certain $\e_0$). Then (by the argument between (2.5) and (2.6)) in $B(0,\frac52)\backslash B(0,\frac 12)$, $E$ is a $C^1$ version of $T$. Therefore in $B(0,\frac52)\backslash B(0,\frac 12)$, $E$ is composed of six $C^1$ faces $E_{ij}$, that are very close to the $T_{ij}$. The $E_{ij}$ meet by three, on four $C^1$ curves $\eta_i,1\le i\le 4$, each $\eta_i$ is very near the half line $[o,y_i)$, and near each $\eta_i$, there exists a tubular neighborhood ${\cal T}_i$ of $\eta_i$, which contains $B([oy_i),r)$ for some $r>0$, in which $E$ is a $C^1$ version of a $\Y$ set. See \cite{DJT} Section 18 for more detail. In total, there is a $C^1$ diffeomorphism $\varphi$, which is very near the identity, such that in $B(0,\frac52)\backslash B(0,\frac 12)$, $E$ coincides with $\varphi(T)$, and $E_{ij}$ corresponds to $\varphi(T_{ij})$, $\eta_i$ corresponds to $\varphi([0,y_i))$.

In particular, since $E$ is very near $T$ in $B(0,\frac52)\backslash B(0,\frac 12)$, $s_{ij}\cap E=\emptyset$, $B_{ij}\cap E=B_{ij}\cap E_{ij}$ is also a one point set, so that locally each $s_{ij}$ links $E_{ij}$, and hence is an element (possibly zero) in $H_1(\R^4\backslash E,\mathbb Z)$, too.

Now we discuss the values in $H_1(\R^4\backslash E,\mathbb Z)$ for these $s_{ij}$.

\begin{lem}
Let $E$ be an Al-minimal set that verifies (6.1). Take all the notations above. Then if
\be \mbox{for all }1\le i<j\le 4, [\vec{s_{ij}}]\ne 0\mbox{ in }H_1(\R^4\backslash E,\mathbb Z),\ee
and 
\be\mbox{ at least 5 of the }[s_{ij}]\mbox{ are mutually different in }H_1(\R^4\backslash E,\mathbb Z),\ee
then $E$ contains at least a point of type other than $\P$ and $\Y$.
\end{lem}
 
\nd We prove it by contradiction. Suppose that there is only $\P$ and $\Y$ points. Then for all $x\in E$, the density $\theta(x)=\lim_{r\to 0}\frac{H^2(B(x,r)\cap E)}{r^2}$ of $E$ at $x$ is either $\frac 32$, either 1. In other words, all singular points in $E$ are of type $\Y$.

Denote by $E_Y$ the set of all the $\Y$ points of $E$. Then $E_Y\cap B(0,2)$ are composed of $C^1$ curves, whose endpoints belong to $\partial B(0,2)$.(c.f.\cite{DJT}, Lemma 18.11, and the $C^1$ regularity around $\Y$ points, see \cite{DEpi} Theorem 1.15 and Lemma 14.6).

The following argument is the same as after Lemma 18.11 in  \cite{DJT} Section 18 (where the reader can get more detail). Here we only sketch the argument.

Since $E$ looks very much like $T$ in $B(0,\frac52)\bs B(0,\frac 12)$, we have $E_Y\cap \partial B(0,2)=\{a_1,a_2,a_3,a_4\}, E_Y\cap\partial B(0,1)=\{b_1,b_2,b_3,b_4\}$, where $b_i$ is the point nearest to $a_i$ among the $b_j,1\le j\le 4$. Then through each $a_i$ passes a curve in $E_Y$, and hence locally $a_i$ lies in the intersection of three half surfaces $E_{ij},j\ne i,1\le j\le 4$. 

But on the sphere $\partial B(0,2)$ we have four $\Y$ points, hence without loss of generality, we can suppose that there is a curve $\gamma_1$ of $E_Y$ that enters the ball $B(0,2)$ by $a_1$ and leaves the ball by $a_2$, and another curve $\gamma_2$ which enters the ball by $a_3$ and leaves it at $a_4$ (see Figure 6-2, where the green curves represent the $\gamma_i,i=1,2$, and we do not know much about the structure of $E$ in $B_{1/2}=B(0,\frac12)$).

          \centerline{\includegraphics[width=0.7\textwidth]{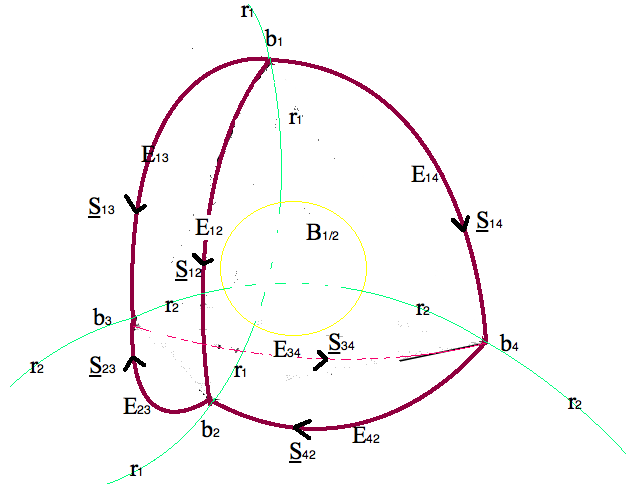}}
          \nopagebreak[4]
          \centerline{6-2} 
          
Near each point $x$ of $\gamma_1$, there exists a $C^1$-ball $B(x,r_x)$ of $x$, in which $E$ is the image of a $\Y$ set under a $C^1$ diffeomorphism. By compactness of the curve $\gamma_1$, there exists a tubular neighborhood $I_1$ of $\gamma_1$, such that $E\cap I_1$ is composed exactly of three surfaces that meet along the curve $\gamma_1$. 

But $\gamma_1$ connects $a_1$ and $a_2$, hence it passes through $b_1$ and $b_2$. But near each $b_i$, the set $E$ is composed of three half surfaces $E_{ij},j\ne i$. Then since $E$ is locally composed of three half surfaces all along $\gamma_1$, these three half surfaces connect $E_{12},E_{13}, E_{14}$ to $E_{21}, E_{2 3}, E_{24}$. Hence we know that $s_{12},s_{13},s_{14}$ are homotopic to $s_{21},s_{23},s_{24}$ (but we do not know which one is homotopic to which one). A similar argument gives also that $
s_{31},s_{32},s_{34}$ are homotopic to $s_{41},s_{42},s_{43}$.  

For the part $\gamma_1$ we have six cases as follows.
\be \begin{array}{ccc}s_{12}\sim s_{21},&s_{13}\sim s_{23},&s_{14}\sim s_{24};\\
s_{12}\sim s_{21},&s_{13}\sim s_{24},&s_{14}\sim s_{23};\\
s_{12}\sim s_{23},&s_{13}\sim s_{21},&s_{14}\sim s_{24};\\
s_{12}\sim s_{23},&s_{13}\sim s_{24},&s_{14}\sim s_{21};\\
s_{12}\sim s_{24},&s_{13}\sim s_{21},&s_{14}\sim s_{23};\\
s_{12}\sim s_{24},&s_{13}\sim s_{23},&s_{14}\sim s_{21}.\end{array}\ee
Note also that automatically $s_{12}\sim s_{21}$, hence the six cases reduce to four (modulo the symmetry between the indices 3 and 4)
\be \begin{array}{l}s_{13}\sim s_{23},s_{14}\sim s_{24};\\
s_{13}\sim s_{24},s_{14}\sim s_{23};\\
s_{12}\sim s_{23}\sim s_{13},s_{14}\sim s_{24};\\
s_{12}\sim s_{23}\sim s_{14},s_{13}\sim s_{24}.\end{array}\ee

Similarly, for the part $\gamma_2$, we have the following four cases
\be \begin{array}{l}s_{31}\sim s_{41},s_{32}\sim s_{42};\\
s_{31}\sim s_{42},s_{32}\sim s_{41};\\
s_{34}\sim s_{41}\sim s_{31},s_{32}\sim s_{42};\\
s_{34}\sim s_{41}\sim s_{32},s_{31}\sim s_{42}.\end{array}\ee

Combine (6.8) and (6.9), we have eight cases
\be \begin{array}{l}1^\circ\ s_{13}\sim s_{23}\sim s_{42}\sim s_{14};\\
2^\circ\ s_{13}\sim s_{23}\sim s_{42}\sim s_{14}\sim s_{43};\\
3^\circ\  s_{13}\sim s_{24}, s_{14}\sim s_{23};\\
4^\circ\ s_{34}\sim s_{41}\sim s_{32},s_{13}\sim s_{24};\\
5^\circ\ s_{13}\sim s_{23}\sim s_{42}\sim s_{14}\sim s_{12};\\
6^\circ\  s_{13}\sim s_{24}, s_{12}\sim s_{23}\sim s_{14};\\
7^\circ\ s_{12}\sim s_{13}\sim s_{23}\sim s_{42}\sim s_{14}\sim s_{43};\\
8^\circ\ s_{13}\sim s_{24},s_{12}\sim s_{14}\sim s_{23}\sim s_{34}.
\end{array}\ee

In particular, at most four of the $[s_{ij}],1\le i<j\le 4$ are different, which contradicts our hypothesis that at least five of the $\{[s_{ij}],1\le i<j\le 4\}$ are different in $H^1(\R^4\backslash E;\mathbb Z)$.\qed

\begin{cor}
Let $E$ be a 2 dimensional reduced Almgren minimal set in $\R^4$ such that (6.1),(6.5) and (6.6) hold. Suppose also that (5.1) holds. Then $E$ is a $\T$ set parallel to $T$.
\end{cor}

\nd By Lemma 6.4, $E$ contains a point $x$ of type other than $\P$ and $\Y$, hence by (5.1) the density $\theta(x)$ of $E$ at $x$ is larger than or equal to $d_T$. Define $\theta(t)=t^{-2}H^2(E\cap B(x,t))$ the density function of $E$ at $x$. By Proposition 5.16 of \cite{DJT}, $\theta(t)$ is non-decreasing on $t$. Then (6.1) and Lemma 16.43 of \cite{DJT} give that $\lim_{t\to\infty}\theta_t=d_T$. But we already know that $\theta(x)=\lim_{t\to 0}\theta(t)\ge d_T$, hence the monotonicity of $\theta$ yields that $\theta(t)=d_T$ for all $t>0$. By Theorem 6.2 of \cite{DJT}, the set $E$ is a minimal cone centered at $x$, with density $d_T$. Thus by (6.1), $E$ is a $T$ centered at $x$.\qed

\medskip

After Corollary 6.11, the rest for us to do is to discuss the case where $E$ is topologically minimal, and no more than 4 of the $[s_{ij}]$ are different. 

First we prove some properties for these $s_{ij}$.

\begin{lem}
1)\be\label{somme} \sum_{j\ne i}[\vec s_{ij}]=0\mbox{ for all }1\le i\le 4.\ee 
2)
For each $i\ne j\ne k$, 
\be\label{zero} [\vec s_{ij}]\ne 0,\ee and
\be\label{voisin}[\vec s_{ij}]\ne[\vec s_{jk}].\ee
\end{lem}

\nd 1) Fix a $1\le i\le 4$. 

We write $\R^4=\R^3\times\R$, where $T\subset\R^3$. 

Recall that $y_i,1\le i\le 4$ are the four $\Y$ points of $T\cap\partial B(0,1)$, $T_{ij}$ is the sector of $T$ passing through the origin and $y_i,y_j$, $x_{ij}$ is the middle point of the great arc passing through $y_i,y_j$, $P_{ij}$ is the plane passing containing $x_{ij}$ and orthogonal to $T_{ij}$, and $s_{ij}=\partial B_{ij}$, where $B_{ij}=B(x_{ij},\frac{1}{10})\cap P_{ij}$.

Denote by $Y_i$ the cone over $Z_i:=\cup_{j\ne i}\widehat{y_ix_{ij}}$, where $\widehat{y_ix_{ij}}$ denotes the great arc connecting  $y_i$ et $x_{ij}$ (see Figure 6-3 of $Y_1\subset\R^3$ below), $C_T$ the convex hull of $Y_i$. Set $C=C_T\times\R$. Since $C$ is a cone, $C\bs T$ is also a cone. Note that $Z_i\subset S^3\cap C$ is a spherical $\Y$ set of dimension 1. We want to show that $\sum_{j\ne i}[\vec s_{ij}]=0$ in $H_1(C\bs T,\Z)$.  

\centerline{\includegraphics[width=0.5\textwidth]{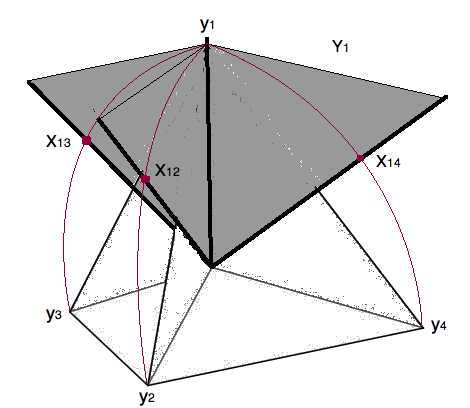}}
\nopagebreak[4]
\centerline{6-3}

Note that $\vec s_{ij}$ is homotopic in $C\bs T$ to its radial projection $\vec s_{ij}'$ on $S^3$ (the orientation of $\vec s_{ij}'$ is induced by $\vec s_{ij}$ on the sphere $S^3$). In fact, denote by $\pi_S$ the radial projection of $\R^4\bs \{0\}$ to $S^3$, then for each $x\in s_{ij}$, the segment $[x,\pi_S(x)]$ belongs to a radial half line, that does not meet any other radial half lines. In particular, since $x\in\R^4\bs T$, where $T$ is a union of radial half lines, $[x,\pi_X(x)]\cap T=\emptyset$. So if we set $f_t(x)=(1-t)x+t\pi_S(x),0\le t\le 1$, then $f_t$ is a homotopy between $\vec s_{ij}$ and $\vec s'_{ij}=\pi_S(\vec s_{ij})$.

Therefore on the sphere, in $C\cap S^3$, the $s_{ij},j\ne i$ are topologically three circles that link respectively the three branches of $Z_i$. Recall that the pair of topological spaces $(C\cap S^3,Z_i)$ is homotopic to $(\R^3,Y)$ where $Y$ is a 1-dimensional $\Y$ set. However in $(\R^3,Y)$, the union of the three oriented circles that link the three branches of $Y$ is the boundary of an oriented manifold with boundary contained in $\R^3\bs Y$. Hence similarly, there exists an oriented manifold with boundary $\Sigma\subset C\cap S^3\bs Z_i$ of dimension 2 such that $\partial \Sigma=\cup_{j\ne i}\vec s'_{ij}$ (see Figure 6-4, where $s_{ij}$ denotes the oriented circle $\vec s_{ij}$, and the orientation of $\Sigma$ is marked by the exterior normal vector $\vec n$). Therefore, after a smooth triangulation under which $\Gamma$ and $s_{ij}$ are all smooth chains, we have $\partial [\Sigma]=\cup_{j\ne i}[\vec s'_{ij}]$. But $\Sigma\subset C\cap S^3\bs Z_i\subset \R^4\bs T$, therefore $\sum_{j\ne i}[\vec s'_{ij}]=0$ in $H_1(\R^4\bs T,\Z)$. Then since $\vec s'_{ij}$ is homotopic to $\vec s_{ij}$,  
\be \sum_{j\ne i}[\vec s_{ij}]=0\mbox{ dans }H_1(\R^4\bs T,\Z).\ee

Now since $E$ is as near as we want to $T$, we can suppose that $\Sigma$ and the $f_t(s_{ij})$ do not touch $E$. Thus we get (\ref{somme}).

\centerline{\includegraphics[width=0.7\textwidth]{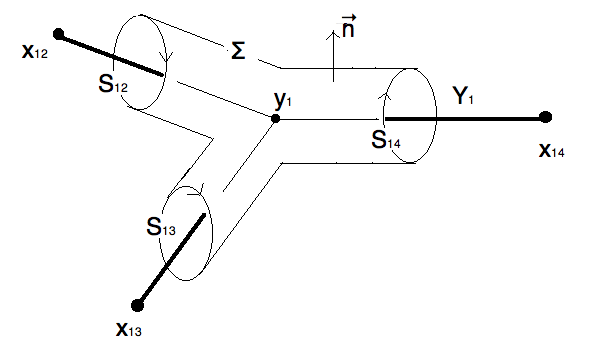}} 
\nopagebreak[4]
\centerline{6-4}

2) Without loss of generality, suppose for example that $i=1,j=2,k=3$. If $[\vec s_{12}]=[\vec s_{23}]$, then by (\ref{somme})
\be [\vec s_{24}]=0.\ee

So we just have to prove (\ref{zero}). 

Suppose for example that $i=2,j=4$. Then (6.17) means that there exists a smooth simplicial 2-chain  $\Gamma$ in $\R^4\bs E$ such that $\partial\Gamma=\vec s_{24}.$ But $E$ is closed, hence there exists a neighborhood $U$ of $\Gamma$ such that $U\cap E=\emptyset$. In particular, $s_{24}\subset U$.
 
Set $D=E\cap \overline B(x_{24},\frac{1}{10})$. Then by the regularity of the minimal set $E$ which is very near $T$, in $B(x_{24},\frac 18)$, $E$ is a piece of very flat surface that is almost a disc. Therefore $D$ is a surface with positive measure.

Set $F=E\bs D$, then $F\bs B(0,2)=E\bs B(0,2)$. We want to show that $F$ is a topological competitor of $E$ with respect to the ball $B(0,2)$. 

So suppose that $\gamma\subset \R^4\bs (B(0,2)\cup E)$ is an oriented circle. We have to show that if $[\gamma]$ is zero in $H_1(\R^4\bs F,\Z)$, then it is zero in $H^1(\R^4\bs E,\Z)$.

Now if $[\gamma]$ is zero in $H_1(\R^4\bs F,\Z)$, then there exists a smooth simplicial 2-chain $\Sigma\subset \R^4\bs F$ such that $\partial \Sigma=\gamma$. By the transversality theorem (c.f. for example \cite{Hir} Chapt 3 Theorem 2.1), we can ask that $\Sigma$ is transversal to $\partial B(x_{24},\frac{1}{10})$.

If $\Sigma\cap B(x_{24},\frac{1}{10})=\emptyset$, then $\Sigma\subset \R^4\bs E$ too, and hence $[\gamma]=0\in H_1(\R^4\bs E,\Z)$. If $\Sigma\cap B(x_{24},\frac{1}{10})\ne\emptyset$, then by the transversality between $\Sigma$ and $\partial B(x_{24},\frac{1}{10})$, and by Proposition 2.36 of \cite{topo}, their intersection is a closed smooth simplicial 1-chain $s\subset\partial B(x_{24},\frac{1}{10})$.

Now we put ourselves in $\overline B_1:=\overline B(x_{24},\frac{1}{10})$. Since $D$ is a very flat topological disc, 
\be H_1(\overline B_1\bs D)=\Z,\ee
whose generator is $[\vec s_{24}]$. As a result, there exists $n\in \Z$ such that $[s]=n [\vec s_{24}]$. Hence there exists a smooth simplicial 1-chain $R\subset \overline B_1\bs D$ such that $\partial R=s-n\vec s_{24}.$ 

Recall that $\Gamma\subset \R^4\bs E$ is such that $\partial \Gamma=\vec s_{24}$. As a result, $\Sigma'=\Sigma\bs \overline B_1+n\Gamma+R$ is a 2-chain satisfying $\partial[\Sigma']=[\gamma]$. Moreover $\Sigma'\subset \R^4\bs E$. Hence $[\gamma]$ is also zero in $H^1(\R^4\bs E,\Z)$. 

Thus, $F=E\bs D$ is a topological competitor of $E$. 

But $D$ is of positive measure, hence
\be H^2(F)<H^2(E),\ee
which contradicts the fact that $E$ is topologically minimal. 

Thus we have proved (\ref{zero}), and hence (\ref{voisin}), and completed the proof of Lemma 6.12.\qed

\bigskip

Now return to our discussion of the case where $E$ is very near a $\T$ set $T$ at scale 1, but contains no point of type other than $\P$ and $\T$. After a discussion (see \cite{XY10} Section 19 for detail) of the 8 cases in (6.10), using Lemma 6.12, the only possibility for $[\vec s_{ij}]$ in $H_1(\R^4\bs E,\Z)$ is: 

\be[\vec s_{13}]=-[\vec s_{24}]=\a,[\vec s_{14}]=-[\vec s_{23}]=\beta,[\vec s_{34}]=\a-\beta, [\vec s_{12}]=-\a-\beta.\ee

Thus we get the following proposition.

\begin{pro}Let $E$ be a reduced topological minimal set of dimension 2 in $\R^4$, that verifies (6.1). Admit all the convention and notation at the beginning of this section, and suppose also that $\gamma_1$ connects $a_1,a_2$, $\gamma_2$ connects $a_3,a_4$. Then if there exists $r>0$ such that $d_{0,3r}(E,T)<\e_0$ (where $\e_0$ the one after (6.3)), but the $s_{ij}$ do not verify (6.20) with respect to $\frac1r E$, then $E$ is a $\T$ set parallel to $T$.
\end{pro}

\section{An example}

In this section, we admit all the notation and convention in Section 6.  We are going to give an example for a set that satisfies (6.20).

Set $w_{ij}=E_{ij}\cap\partial B(0,1)$, (See Figure 7-1, where $\underline w_{ij}$ denote $\vec w_{ij}$), then the $w_{ij}$ are $C^1$ curves. Denote also by $\vec w_{ij}$ the oriented curve from $b_i$ to $b_j$.

\centerline{\includegraphics[width=0.7\textwidth]{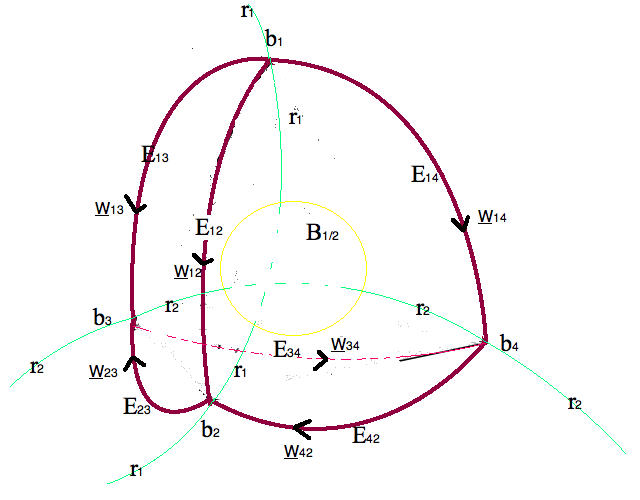}}
\nopagebreak[4]
\centerline{7-1}

Now suppose that $E_Y=\gamma_1\cup\gamma_2$. In other words, all points in $E$ are $\P$ points, except for the two curves. (For the case where $E_Y\ne\gamma_1\cup\gamma_2$, we know that $E_Y\bs(\gamma_1\cup\gamma_2)$ is a union of closed curves, because the only endpoints of $E_Y\cap \partial B(0,1)$ are $\{b_i\}_{1\le i\le 4}$.  This is thus a more complicated case.) 

\begin{lem} $\gamma_1\cup \gamma_2\cup w_{12}\cup w_{34}$ is the boundary of a $C^1$ surface $S_0\subset E$, and $S_0$ contains only points of type $\P$.
\end{lem}

\nd By the $C^1$ regularity of minimal sets,the part of $E$ in $B(0,1)$ is composed of $C^1$ manifolds $S_1,S_2,\cdots$ whose boundaries are unions of curves in the set $Bd=\{w_{ij},\gamma_1,\gamma_2\}$. Thus there exists $k\in\N$ such that $w_{12}$ is part of the boundary of the manifold $S_k$. But $\partial S_k$ is a union of several closed curves, while $w_{12}$ is not closed. Hence there exists a curve $\gamma\in Bd$ that touches $w_{12}$ and such that $\gamma$ is also part of $\partial S_k$. If one of the $w_{1i}$ (resp. one of the $w_{2j}$) is part of $\partial S_k$ and touches $w_{12}$,  we have $[\vec s_{12}]=[\vec s_{i1}]$ (with orientation) (resp. $[\vec s_{12}]=[\vec s_{2i}]$), which contradicts (\ref{voisin}).

Hence the only possibility for $\gamma$ is $\gamma_1$.  This means, the union of $w_{12}$ and $\gamma_1$ is part of the boundary of a manifold $S_k$, and except for $w_{34}$, the boundary of $S_k$ contains no other $w_{ij}$. A similar argument gives also that the union of $w_{34}$ and $\gamma_2$ is part of the boundary of a manifold $S_l$, and the boundary of $S_l$ contains no other $w_{ij}$, except perhaps for $w_{12}$.

Thus, either the union of these four curves $w_{12}, w_{34},\gamma_1,\gamma_2$ is the boundary of a surface $S_k$, either the union of $w_{12}$ and $\gamma_1$ and the union of $w_{34}$ and $\gamma_2$ are the boundaries of two surfaces $S_k$ and $S_l$. In any case, the union of the four curves is the boundary of a $C^1$ surface $S_0\subset E$, which is not necessarily connected.\qed

After Lemma 7.1, if we take away the surface $S_0$ from $E$, then $E\bs S$ is composed of a union of $C^1$ surfaces, whose boundaries are unions of curves belonging to $Bd'=\{w_{13},w_{14},w_{23},w_{24},\gamma_1,\gamma_2\}$. Still by the same argument above, there are two surfaces $S_1,S_2$, with $\partial S_1=w_{13}\cup \gamma_2\cup w_{24}\cup \gamma_1$, and the other one $\partial S_2=w_{23}\cup \gamma_2\cup w_{14}\cup \gamma_1$. Moreover $S_1\cup S_2$ is also a connected topological manifold, for which we can define local orientation, even near $\partial S_1$ and $\partial S_2$.

\centerline{\includegraphics[width=0.8\textwidth]{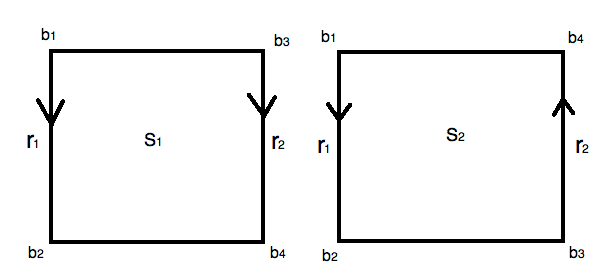}}
\nopagebreak[4]
\centerline{7-2}

Thus topologically, the boundaries of the two surfaces $S_1,S_2$ are like the boundaries of two squares, one with the four vertices $b_1,b_3,b_4,b_2$ (where we write them in the order of adjacent), the other with four vertices  $b_1,b_4,b_3,b_2$. Moreover, we have to glue the two $\overrightarrow{b_1b_2}$ in $S_1$ and $S_2$ together, and the same for $\overrightarrow{b_3b_4}$. Note that these two gluings are of different directions (See Figure 7-2).

Notice also that after the gluing, $S_1\cup S_2$ cannot be orientable.

\begin{rem}Since $S_1\cup S_2$ is not orientable, $[\vec s_{13}],[\vec s_{14}],[\vec s_{24}],[\vec s_{23}]$ are all of order 2 in $H_1(\R^4\bs E,\Z)$. 

In fact for a connected surface $S$, the non-orientability means that for each point $x\in S$ we can find a path $\gamma:[0,1]\to S$ such that $\gamma(0)=\gamma(1)=x$, and if we denote by $n(t)=x(t)\wedge y(t)\in \wedge_2 N_{\gamma(t)}S$ a continuous unit normal 2-vector field on $\gamma$, where $x(t),y(t)\in N_{\gamma(t)}S$ are unit normal vector fields, and $n,x,y$ are continuous with respect to $t$, then $n(0)=-n(1).$ Note that $n(t)$ can also represent the oriented plane in $\R^4$. Define, for each $r>0$ $s_r(t):T=\R/\Z\to P_t=P(x(t)\wedge y(t)), \theta\mapsto r[\cos (2\pi \theta) x(t)+\sin (2\pi\theta) y(t)]$. The the image of $s_r(0)$ and $s_r(1)$ are the same circle, but with different orientations: $s_r(0)(t)=s_r(1)(-t)$.

Set $Q_t=\gamma(t)+P(t)$. Fix $r>0$ sufficiently small, such that for each $t\in[0,1]$, $B(\gamma(t),r)\cap Q_t\cap S=\{\gamma(t)\}.$ 

Set $G:T\times[0,1]\to\R^4, G(\theta,t)=s_r(t)(\theta)+\gamma(t)$. This is a continuous map, with $G(T\times\{0\})=s_r(0)$ and $G(T\times\{1\})=s_r(1)=-s_r(0)$. As a result, the oriented circle $s_r(0)$ is homotopic to $-s_r(0)$, and is hence of order 2.

Now for each $s\in\{[\vec s_{13}],[\vec s_{14}],[\vec s_{24}],[\vec s_{23}]\}$, we can first find a circle $s'$ homotopic to $s$, such that there exists $x,\gamma$ as before, and that there exists $R>0$ such that $s_R(0)=s'$. We can find $r>0$ as above, then $s_r(0)$ is homotopic to $s'$, and hence $s$. Therefore $[s]=[s_r(0)]$ is of order 2.
\end{rem}

We will construct a set $E\subset \R^4$, with all the above properties. That is, in $B(0,1)$, the set $E$ is the union of $S_0$ and $S_1\cup S_2$ as above, $S_1\cup S_2$ is a non-orientable topological manifold; $S_0$ has two connected components, that meet $S_1\cup S_2$ at $\gamma_1$ and $\gamma_2$ respectively. Outside the ball $B(0,1)$, $E$ is a $C^1$ version of $T$, and it looks like $T$ at infinity. Moreover, $H_1(\R^4\bs E)$ and the $[\vec s_{ij}]$ satisfy (6.20).

Take two copies of squares (see Figure 7-2), one with vertices (written in the order of the sense of clockwise) $b_1,b_3,b_4,b_2$, the other with vertices $b_1,b_4,b_3,b_2$. We glue the two sides $\overrightarrow{b_3b_4}$ in $S_1$ and $S_2$ together, and we do the same for $\overrightarrow{b_1b_2}$. Thus we get a M\"obius band in $\R^3$.  (See Figure 7-3.)

\centerline{\includegraphics[width=0.7\textwidth]{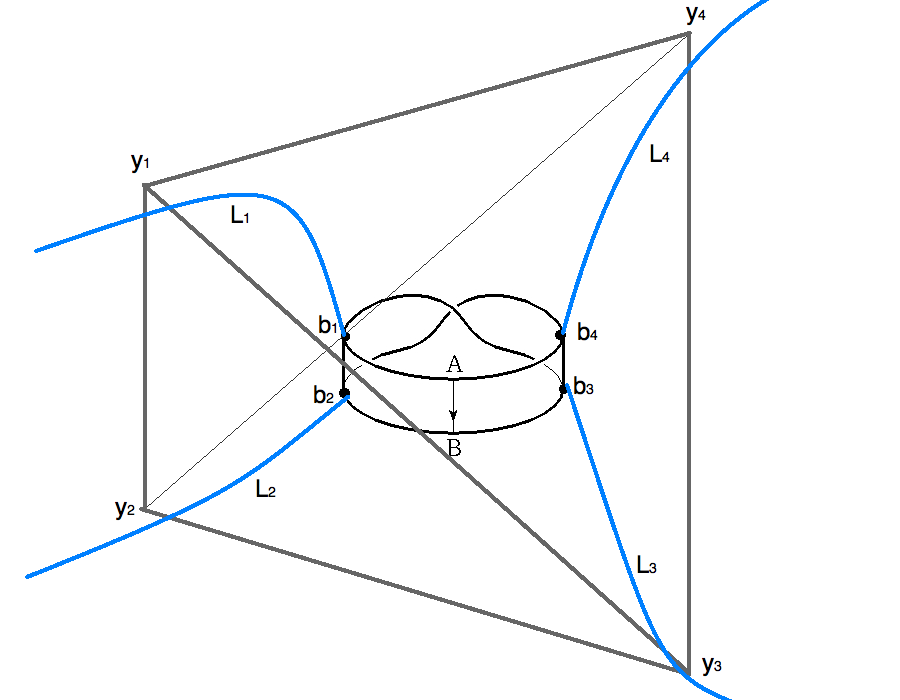}}
\nopagebreak[4]
\centerline{7-3}

Next, take a very big regular tetrahedron centered at the origin (with vertices $y_i,1\le i\le 4$) which contains the M\"obius band constructed before.  For each $i$, take a smooth curve $L_i$ issued from $b_i$ and goes to infinity, such that $L_i$ tends to $[0,y_i)$. (See Figure 7-3).

Then take, for each $1\le i\ne j\le 4$, a $C^1$ surface $E_{ij}$, homeomorphic to $\R^2$, whose boundary is $L_i\cup L_j\cup [b_ib_j]$. Note that all $E_{ij}$ go to infinity, hence in $\R^3$, $E_{23}$ and $E_{14}$, or $E_{13}$ and $E_{24}$ has to meet each other. Thus we move to $\R^4$ to avoid this. 

Thus we get a set that looks like a $T$ at infinity, and at first glance, we cannot say easily that a set with such a topology cannot be topologically minimal.

\renewcommand\refname{References}
\bibliographystyle{plain}
\bibliography{reference}

\end{document}